\input amstex
\input epsf
\input cyracc.def
\font\twcyr=wncyr10 scaled \magstep1

\documentstyle{amsppt}\nologo\footline={}
\hsize450pt

\def\B{\mathop{\text{\rm B}}}
\def\E{\mathop{\text{\rm E}}}

\def\S{\mathop{\text{\rm S}}}
\def\T{\mathop{\text{\rm T}}}
\def\d{\mathop{\text{\rm d}}}
\def\arccosh{\mathop{\text{\rm arccosh}}}
\def\per{\mathop{\text{\rm per}}}
\def\sgn{\mathop{\text{\rm sgn}}}
\def\ta{\mathop{\text{\rm ta}}}
\def\vol{\mathop{\text{\rm vol}}}
\def\Lin{\mathop{\text{\rm Lin}}}

\def\lob{\mathop{\text{\twcyr l}}}

\topmatter\title Seidel's conjectures in hyperbolic 3-space\endtitle
\author Omar Chavez Cussy, Carlos H.~Grossi\endauthor
\thanks First author supported by CAPES.\endthanks
\address Departamento de Matem\'atica, ICMC, Universidade de S\~ao Paulo,
Caixa Postal 668,\newline13560-970--S\~ao Carlos--SP, Brasil\endaddress
\email omarchavez\@usp.br\endemail
\address Departamento de Matem\'atica, ICMC, Universidade de S\~ao Paulo,
Caixa Postal 668,\newline13560-970--S\~ao Carlos--SP, Brasil\endaddress
\email grossi\@icmc.usp.br\endemail
\abstract We prove, in the case of hyperbolic $3$-space, a couple of conjectures raised by J.~J.~Seidel in
{\it On~the volume of a hyperbolic simplex,} Stud.~Sci.~Math.~Hung.~{\bf 21}, 243–-249, 1986. These
conjectures concern expressing the volume of an ideal hyperbolic tetrahedron as a monotonic function of
algebraic maps. More precisely, Seidel's first conjecture states that the volume of an ideal tetrahedron in
hyperbolic 3-space is determined by (the permanent and the determinant of) the doubly stochastic Gram matrix
$G$ of its vertices; Seidel's fourth conjecture claims that the mentioned volume is a monotonic function of
both the permanent and the determinant of $G$.
\endabstract\endtopmatter\document

\vskip-4pt

\centerline{\bf1.~Introduction}

\medskip

It is a typical phenomenon in hyperbolic geometry\footnote{As well as in several other {\it classic\/}
geometries (see [AGr2]).} that explicit formulae for calculating even simple geometric invariants involve
transcendental functions. These functions can be quite sophisticated and such is usually the case when it
comes to volume-related problems which are our main concern in this paper. A way to deal with this
difficulty is to express the geometric invariants in question as {\it monotonic functions\/} of {\it
algebraic maps\/}: in a certain sense, this allows us to `replace' complicated geometric invariants by much
simpler ones.

As a toy example, let us consider the distance function in the projective model of hyperbolic $n$-space.
Take an $(n+1)$-dimensional $\Bbb R$-linear space $V$ equipped with a bilinear symmetric form of signature
$--\dots-+$. The hyperbolic $n$-space is the open $n$-ball of positive points
$\Bbb H^n:=\big\{\pmb p\in\Bbb PV\mid\langle p,p\rangle>0\big\}$.
(We denote respectively by $\pmb p$ and $p$ a point $\pmb p\in\Bbb PV$ and a representative
$p\in V$.) The {\it ideal\/} boundary of~$\Bbb H^n$, also known as the {\it absolute,} is the
$(n-1)$-sphere $\partial\Bbb H^n:=\big\{\pmb p\in\Bbb PV\mid\langle p,p\rangle=0\big\}$. Hyperbolic
$n$-space is endowed with the distance function
$$d(\pmb p,\pmb q):=\arccosh\displaystyle\sqrt{\frac{\langle p,q\rangle\langle q,p\rangle}{\langle
p,p\rangle\langle q,q\rangle}},\quad\pmb p,\pmb q\in\Bbb H^n.$$
Clearly, distance is a monotonic function of the {\it tance\/}
$\ta(\pmb p,\pmb q):=\displaystyle\frac{\langle p,q\rangle\langle q,p\rangle}{\langle
p,p\rangle\langle q,q\rangle}$.
Due to its algebraic nature, it is usually much simpler to use the tance instead of the distance in
applications (see, for instance, [AGG], [AGr2], [Ana1], [Ana2]). J.~J.~Seidel's conjectures [Sei] concern
applying a similar ideia to the case of the volume of an ideal simplex in $\Bbb H^n$.

A (labelled) {\it ideal simplex\/} in $\Bbb H^n$ is an $(n+1)$-tuple $(\pmb v_1,\dots,\pmb v_{n+1})$ of
ideal points $\pmb v_i\in\partial\Bbb H^n$ called the {\it vertices\/} of the ideal simplex. The volume of
an ideal simplex is the hyperbolic volume of the convex hull of the points $\pmb v_1,\dots,\pmb v_{n+1}$.
In dimension $3$, for example, the volume of an ideal simplex --- an ideal {\it tetrahedron\/} --- can be
expressed (see Section 5) in terms of its dihedral angles and the Lobachevsky function
$\lob:\Bbb R\to\Bbb R$,
$$\lob(\theta):=-\int_0^{\theta}{\log|2\sin t|\d\!t}.\leqno{(1.1)}$$

Let $S:=(\pmb v_1,\dots,\pmb v_{n+1})$ be an ideal simplex in $\Bbb H^n$. Choosing representatives
$v_i\in V$ we obtain a~{\it Gram matrix\/} $G$ of the vertices of $S$, where $G:=[g_{ij}]$,
$g_{ij}:=\langle v_i,v_j\rangle$. Among all the Gram matrices of the vertices of $S$, there is a single
one\footnote{Some particular (very) degenerate ideal simplices do not admit a doubly stochastic Gram matrix
(see Remark 3.2.3).} that is {\it doubly stochastic\/} (a square matrix is called doubly stochastic if all
its entries are non-negative and the sum of entries in every row and in every column equals $1$).
Let~$G_{ds}$ stand for the doubly stochastic Gram matrix of the vertices of $S$.

\smallskip

J.~J.~Seidel's first conjecture is that the volume of $S$ should be determined by (some natural algebraic
functions of) the entries of $G_{ds}$. The fourth conjecture states that the volume of $S$ is a monotonic
function of both the determinant and the permanent\footnote{The permanent of an $m\times m$ matrix
$A=[a_{ij}]$ is defined by the expression
$\per A:=\sum_{\sigma\in S_m}a_{1\sigma(1)}a_{2\sigma(2)}\dots a_{m\sigma(m)}$, where $S_m$ stands for the
symmetric $m$-group.} of $G_{ds}$. We prove these conjectures in the case of hyperbolic $3$-space. Two
observations are in order here. First, we establish in Theorem 4.1 a stronger version of the first
conjecture by showing that permanent and determinant actually serve as coordinates of the space of ideal
tetrahedra modulo isometries (we also provide in Theorem 5.1 an explicit formula for the volume of an
ideal tetrahedron in terms of the permanent and the determinant of $G_{ds}$). Secondly, Seidel's original
fourth conjecture says that the volume is a decreasing function of the permanent, while it is in fact an
increasing one (see Theorem 6.10). Trying to understand what may have led Seidel to believe that the volume
is decreasing in the permanent, we proved his third conjecture (see Theorem 6.11 and the related Proposition
6.12) which is a little bit technical and, in the particular low-dimensional case we consider, quite simple
to show.

\smallskip

It should be emphasized that taking the {\it doubly stochastic\/} Gram matrix of {\it the vertices\/} of $S$
seems to be crucial\footnote{Well, it is possible to take $DG_{ds}D$ in place of $G_{ds}$, where $D$ is
a (fixed) diagonal matrix whose diagonal entries are non-null constants, as this only multiplies the
determinant and permanent of every $G_{ds}$ by the same positive number $(\det D)^2$.} to Seidel's
conjectures. In fact, N.~V.~Abrosimov [Abr] proved that the variant of the first conjecture where one
considers, instead of $G_{ds}$, a normalized Gram matrix of the points {\it polar\/} to the faces of $S$
(polar points are discussed in Sections 2.1 and 2.2), actually does not hold in full generality: in this
case, the volume of an ideal tetrahedron can be expressed as a function of determinant and permanent if and
only if the tetrahedron is either acute-angled or obtuse-angled [Abr, Theorem 3 and Example 1]. Moreover,
when we tried to consider Gram matrices of the vertices of $S$ which were not doubly stochastic (there are
other apparently `natural' choices of representatives), the fourth conjecture turned out to be false. 

\smallskip

Seidel cites some results about doubly stochastic Gram matrices of vertices of simplices that are intended
to justify this particular choice of Gram matrices in his conjectures [Sei, Facts 1, 4, 6]. Some of~these
facts involve relationships between the exterior and symmetric algebras of a linear space. Curiously,
determinant and permanent are particular cases of what Littlewood-Richardson called the {\it immanants\/} of
a matrix [LiR]. Immanants are closely related to Schur functors which, in their turn, seem to play an
important role in hyperbolic geometry as well as in other classic geometries. (Exterior and symmetric powers
are the Schur functors corresponding respectively to the determinant and permanent.) For~instance,
calculating the dihedral angles of an ideal tetrahedron in terms of the doubly stochastic Gram matrix of
its vertices involves the use of the exterior power functor and of the Hodge star operator (see Section
2.2 and Proposition 2.2.4).

\smallskip

In~Section 7 we briefly discuss the relationship between immanants and Schur functors; it is reasonable to
expect that other immanants (besides the determinant and permanent) will be involved in possible
generalizations of J.~J.~Seidel's conjectures to higher dimensions and/or to more general polyhedra (not
necessarily ideal ones).

\newpage

\centerline{\bf2.~Preliminaries}

\medskip

{\bf2.1.~Hyperbolic $3$-space.}~The approach to hyperbolic geometry in this section follows [AGr2].

\smallskip

Let $V$ be an $\Bbb R$-linear space equipped with a bilinear symmetric form
$\langle-,-\rangle:V\times V\to\Bbb R$ of signature $---+$. This form divides the projective $3$-space
$\Bbb PV$ into {\it positive,} {\it null\/} (or {\it ideal\/}), and {\it negative\/} points:
$$\B V:=\big\{\pmb p\in\Bbb PV\mid\langle p,p\rangle>0\big\},\quad
\S V:=\big\{\pmb p\in\Bbb PV\mid\langle p,p\rangle=0\big\},\quad
\E V:=\big\{\pmb p\in\Bbb PV\mid\langle p,p\rangle<0\big\}.$$
(As stated in the introduction, we use the notation $\pmb p$ and $p$ respectively for a point
$\pmb p\in\Bbb PV$ and a representative $p\in V$.) Note that $\B V$ is an open $3$-ball whose boundary is
the $2$-sphere $\S V$, called the {\it absolute.}

Let $\pmb p\in\Bbb PV$ be a non-ideal point. We have the well-known natural identification
$${\T}_{\pmb p}\Bbb PV=\Lin(\Bbb Rp,p^\perp)\leqno{(2.1.1)}$$
of the tangent space to $\Bbb PV$ at $\pmb p$ and the $\Bbb R$-linear space of linear maps from the subspace
generated by $p$ to its orthogonal $p^\perp$ with respect to $\langle-,-\rangle$. Given tangent vectors
$\varphi_1,\varphi_2\in\Lin(\Bbb Rp,p^\perp)$ at a point $\pmb p\in\Bbb PV$, we define
$$\langle\varphi_1,\varphi_2\rangle:=-\frac{\big\langle\varphi_1(p),\varphi_2(p)\big\rangle}{\langle
p,p\rangle}.\leqno{(2.1.2)}$$
This endows $\B V$ with a Riemannian metric known as the {\it hyperbolic\/} metric and $\E V$ with a
Lorentzian metric sometimes called the {\it de Sitter\/} metric. We call $\B V$ the {\it hyperbolic\/
$3$-space.}

\smallskip

Any $2$-dimensional linear subspace $W\leqslant V$ of signature $-+$ gives rise to a geodesic
$\Bbb PW\cap\B V$ in $\B V$ and every geodesic in $\B V$ is obtained in this way. (An analogous statement
holds for de Sitter space.) Any $3$-dimensional linear subspace $W\leqslant V$ of signature $--+$ gives rise
to a totally geodesic surface $\Bbb PW\cap\B V$ in $\B V$ known as a {\it plane.} There is a simple duality
between planes in $\B V$ and points in $\E V$: the negative point $\Bbb PW^\perp\in\E V$ corresponds to the
plane $\Bbb PW\cap\B V$, where $W\leqslant V$ is a $3$-dimensional linear subspace of signature $--+$. In
other words, the Lorentzian manifold $\E V$ is simply the space of all planes in $\B V$. The point
$\Bbb PW^\perp\in\E V$ is called the {\it polar\/} point to the plane $\Bbb PW$.

\smallskip

In what follows, we will denote the hyperbolic $3$-space $\B V$ by $\Bbb H^3$. The projective space
$\Bbb PV$ will be referred to as the {\it extended\/} hyperbolic $3$-space.

\bigskip

{\bf2.2.~Angle.}~Let us apply the duality between points and planes in the extended hyperbolic space in
order to find an expression for the angle between two planes. This expression will be used later to find the
dihedral angles of an ideal tetrahedron in terms of the entries of a certain Gram matrix of the tetrahedron
(see Remark 3.2.4). First, we remind the reader a few basic facts about the Hodge star operator. The
results in this subsection are related to those in [Moh, Section 2.2].

\smallskip

Let $U$ stand for an $\Bbb R$-linear space equipped with a non-degenerate bilinear symmetric form
$\langle-,-\rangle:U\times U\to\Bbb R$ of arbitrary signature $(n,m)$, where $N:=\dim U=n+m$ and $n$
denotes the negative part of the signature. Let $\sigma=1$ if $n=0\!\!\mod2$ and $\sigma=-1$ otherwise.

The $k$-th exterior power $\bigwedge^kU$, $1\leqslant k\leqslant N$, is equipped with the bilinear
symmetric form defined by
$$\langle v_1\wedge\dots\wedge v_k,w_1\wedge\dots\wedge w_k\rangle:=\det[g_{ij}],\quad
g_{ij}:=\langle v_i,w_j\rangle.$$
This form on $\bigwedge^kU$ is non-degenerate. Indeed, let
$(b_1,\dots,b_{N})$ be an orthonormal basis for $U$, that is, $\langle b_i,b_i\rangle=\sigma_{i}$ and
$\langle b_i,b_j\rangle=0$ for $i\ne j$, where $\sigma_i=\pm1$. Then
$(b_{i_1}\wedge\dots\wedge b_{i_k}\mid 1\leqslant i_1<\dots<i_k\leqslant N)$ is an orthonormal basis for
$\bigwedge^kU$.

\smallskip

Fix $\omega\in\bigwedge^NU$ such that $\langle\omega,\omega\rangle=\sigma$. The Hodge star operator is the
$\Bbb R$-linear map $*:\bigwedge^kU\to\bigwedge^{N-k}U$, $b\mapsto*b$, defined by requiring that
$a\wedge*b=\langle a,b\rangle\omega$ for every $a\in\bigwedge^kU$. Clearly, $*\omega=\sigma$ since
$\omega\wedge*\omega=\langle\omega,\omega\rangle\omega=\sigma\omega$.

\smallskip

The following proposition is an adaptation of [Huy, Proposition 1.2.20] to the case of a non-degenerate form
of arbitrary signature. The properties of the Hodge star operator in the propositon seem to be well-known
but we could not find a direct reference. The proofs are elementary and provided only for the sake of
completeness.

\medskip

{\bf2.2.1.~Proposition.}~{\it The Hodge star operator is injective. It satisfies the following
identities\/{\rm:}

\smallskip

$\bullet$ $a\wedge*b=b\wedge*a$ for every\/ $a,b\in\bigwedge^kU$

\smallskip

$\bullet$ $\langle a,b\rangle=\sigma\cdot*(a\wedge*b)=\sigma\cdot*(b\wedge*a)$ for every\/
$a,b\in\bigwedge^kU$

\smallskip

$\bullet$ $\langle a,*b\rangle=(-1)^{k(N-k)}\langle *a,b\rangle$ for every\/ $a\in\bigwedge^kU$ and\/
$b\in\bigwedge^{N-k}U$

\smallskip

$\bullet$ $*(*a)=(-1)^{k(N-k)}\sigma\,a$ for every $a\in\bigwedge^kU$}

\medskip

{\bf Proof.}~If $*b=0$, then $a\wedge*b=\langle a,b\rangle\omega=0$ for every $a\in\bigwedge^kU$ which
implies $b=0$ because the form on $\bigwedge^kU$ is non-degenerate.

The first identity is obvious. The second follows from
$*(a\wedge*b)=\langle a,b\rangle\cdot*w=\sigma\langle a,b\rangle$. The third follows from the second:
$$\langle a,*b\rangle=\sigma\cdot*(*b\wedge*a)=(-1)^{k(N-k)}\sigma\cdot*(*a\wedge*b)=(-1)^{k(N-k)}
\langle*a,b\rangle.$$
As for the last equality, take an orthonormal basis $(b_1,\dots,b_N)$ in $U$, that is,
$\langle b_i,b_i\rangle=\sigma_i=\pm1$ and $\langle b_i,b_j\rangle=0$ for $i\ne j$. We can assume that
$\omega=b_1\wedge\dots\wedge b_N$. Fix $i_1<i_2<\dots<i_k$ and let $j_1<j_2<\dots<j_{N-k}$ be the indices
complementary to $i_1<i_2<\dots<i_k$. Let us show that
$$*(b_{i_1}\wedge\dots\wedge b_{i_k})=\sgn(h)\cdot\sigma_{i_1}\cdots\sigma_{i_k}\cdot
b_{j_1}\wedge\dots\wedge b_{j_{N-k}},$$
where $h$ is the permutation $(i_1,i_2,\dots,i_k,j_1,j_2,\dots,j_{N-k})\mapsto(1,2,\dots,N)$. Indeed,
$$(b_{l_1}\wedge\dots\wedge b_{l_k})\wedge\big(\sgn(h)\cdot\sigma_{i_1}\cdots\sigma_{i_k}\cdot
b_{j_1}\wedge\dots\wedge b_{j_{N-k}}\big)=\cases0,\text{ if } (l_1,\dots,l_k)\ne(i_1,\dots,i_k)\\
\sigma_{i_1}\dots\sigma_{i_k}\cdot\omega,\text{ if } (l_1,\dots,l_k)=(i_1,\dots,i_k)\endcases=$$
$$=\langle b_{l_1}\wedge\dots\wedge b_{l_k},b_{i_1}\wedge\dots b_{i_k}\rangle\omega$$
for every $l_1<\dots<l_k$. It remains to observe that
$$*\big(*(b_{i_1}\wedge\dots\wedge
b_{i_k})\big)=\sgn(h)\cdot\sigma_{i_1}\dots\sigma_{i_k}\cdot*(b_{j_1}\wedge\dots\wedge
b_{j_{N-k}})=(-1)^{k(N-k)}\sigma\cdot b_{i_1}\wedge\dots\wedge b_{i_k},$$
where $(-1)^{k(N-k)}=\sgn(h)\cdot\sgn(h')$ and $h'$ stands for the permutation
$(j_{1},\dots,j_{N-k},i_1,\dots,i_k)\mapsto\break(1,2,\dots,N)$.
\hfill$_\blacksquare$

\medskip

From now on, we are back to the particular case of the $\Bbb R$-linear space $V$ of signature $(3,1)$.

\medskip

{\bf2.2.2.~Remark.~{\rm[Moh, Remark 1, Section 2.2]}.} {Let\/ $\pmb v_1,\pmb v_2,\pmb v_3\in\S V$ be
pairwise distinct ideal points and let\/ $P:=\Bbb PW$, $W:=\Bbb Rv_1+\Bbb Rv_2+\Bbb Rv_3$, be the plane
generated by\/ $\pmb v_1,\pmb v_2,\pmb v_3$. Then\/ $\pmb u$ is the polar point to\/ $P$, where\/
$u:=*(v_1\wedge v_2\wedge v_3)$.}

\medskip

{\bf Proof.}~By Proposition 2.2.1,
$\big\langle v_i,*(v_1\wedge v_2\wedge v_3)\big\rangle\omega=v_i\wedge*\big(*(v_1\wedge v_2\wedge
v_3)\big)=v_i\wedge v_1\wedge v_2\wedge v_3=0$ for $i=1,2,3$.
\hfill$_\blacksquare$

\medskip

{\bf2.2.3.~Lemma.}~{\it Let\/ $P$ be a plane with polar point\/ $\pmb u\in\E V$ and let\/ $\pmb p\in P$.
Let\/ $\pmb v\in\B V\cup\S V$ be a point that does not belong to\/ $P$. The tangent vector\/
$n:=\langle-,p\rangle u$, $x\mapsto\langle x,p\rangle u$ {\rm(}see {\rm(2.1.1))}, points towards the
half-space of\/ $\B V\cup\S V$ determined by\/ $P$ and containing\/ $v$ if and only if\/
$\langle u,v\rangle\langle v,p\rangle<0$.}

\medskip

{\bf Proof.} By [AGG, Lemma 4.2.2], the tangent vector $n$ is normal to $P$ at $\pmb p$. It follows from
[AGr2, Lemma 5.2] that $\langle-,p\rangle\displaystyle\frac{\pi[p]v}{\langle v,p\rangle}$ is tangent to the
oriented segment of geodesic from $\pmb p$ to $\pmb v$, where
$\pi[p]v:=v-\displaystyle\frac{\langle v,p\rangle}{\langle p,p\rangle}p$. Then $n$ points towards the
half-space determined by $P$ containing $\pmb v$ if and only if
$$\bigg\langle\langle-,p\rangle u,\langle-,p\rangle\displaystyle\frac{\pi[p]v}{\langle v,p\rangle}
\bigg\rangle=-\frac{\langle p,p\rangle\langle u,v\rangle}{\langle v,p\rangle}>0,$$
where the above product is taken with respect to the Riemannian metric introduced in (2.1.2).
\hfill$_\blacksquare$

\medskip

Let $\pmb v_1,\pmb v_2,\pmb v_3,\pmb v_4\in\S V$ be pairwise distinct ideal points. Let $P_1$ and $P_2$ be
the planes respectively generated by $\pmb v_1,\pmb v_2,\pmb v_3$ and $\pmb v_1,\pmb v_2,\pmb v_4$. The
geodesic $\gamma$ joining $\pmb v_1,\pmb v_2$ is common to $P_1$ and $P_2$. Our~intention is to measure the
(non-oriented) angle in $[0,\pi]$ between the half-plane in $P_1$ which contains $\pmb v_3$ and is
determined by $\gamma$ and the half-plane in $P_2$ which contains $\pmb v_4$ and is determined by $\gamma$.
This angle is called the {\it dihedral angle\/} between $P_1$ and $P_2$ at the edge $\gamma$.

\medskip

{\bf2.2.4.~Proposition.}~{\it Let\/ $\pmb v_1,\pmb v_2,\pmb v_3,\pmb v_4\in\S V$ be pairwise distinct
ideal points and let\/ $g_{ij}:=\langle v_i,v_j\rangle$. There exist representatives\/ $v_i\in V$ so that\/
$g_{ij}>0$ for $i\ne j$. Moreover, with such a choice of representatives, the dihedral angle\/
$\angle(P_1,P_2)$ between the plane\/ $P_1$ generated by\/ $\pmb v_1,\pmb v_2,\pmb v_3$ and the plane\/
$P_2$ generated by\/ $\pmb v_1,\pmb v_2,\pmb v_4$ is
$$\angle(P_1,P_2)=\arccos
\frac{g_{13}g_{24}+g_{14}g_{23}-g_{12}g_{34}}{2(g_{23}g_{31}g_{24}g_{41})^{\frac12}}.$$}

{\bf Proof.}~Clearly, we can assume that $g_{12},g_{23},g_{34}>0$. Since $\pmb v_1,\pmb v_2,\pmb v_3$ are
ideal and pairwise distinct, the form on the subspace generated by $v_1,v_2,v_3$ has
signature $--+$; the same holds for the subspaces generated by $v_1,v_2,v_4$ and by $v_2,v_3,v_4$. Hence,
the Gram matrices of $v_1,v_2,v_3$, of $v_1,v_2,v_4$, and of~$v_2,v_3,v_4$ have positive determinant. This
implies that $g_{ij}>0$ for $i\ne j$.

Assume that $\pmb v_4$ does not belong to $P_1$. Let $\gamma$ be the geodesic joining $\pmb v_1,\pmb v_2$.
Then $p:=v_1+v_2$ gives rise to a point $\pmb p\in\gamma$ since $\langle v_1+v_2,v_1+v_2\rangle=2g_{12}>0$.
We define
$$\omega:=\frac1{\sqrt{-\det G}}\;v_1\wedge v_2\wedge v_3\wedge v_4\in{\bigwedge}^4V,$$
where $G:=[g_{ij}]$ stands for the Gram matrix of $v_1,v_2,v_3,v_4$ (its determinant is negative due to the
signature $---+$ of the form on $V$). Note that $\langle\omega,\omega\rangle=-1$. Let
$u_1:=*(v_1\wedge v_2\wedge v_3)$ and $u_2:=*(v_1\wedge v_2\wedge v_4)$. By Remark 2.2.2, $\pmb u_1$ and
$\pmb u_2$ are respectively the polar points to $P_1$ and $P_2$.

\medskip
\noindent
$\epsfxsize=3cm\vcenter{\hbox{\epsfbox{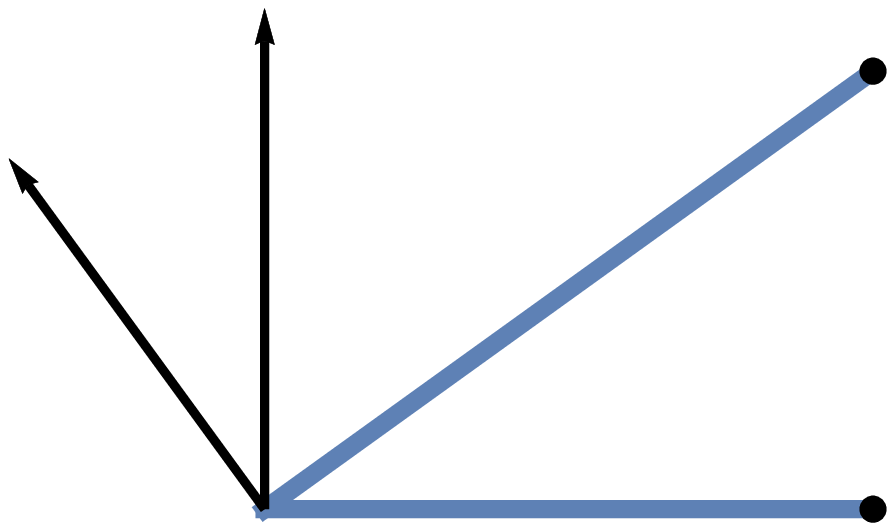}}}$

\leftskip11pt
$\pmb p$

\leftskip69pt\vskip-50pt
$\pmb v_4$

\vskip26pt
$\pmb v_3$

\leftskip39pt\vskip-10pt
$P_1$

\leftskip35pt\vskip-52pt
$P_2$

\leftskip19pt\vskip-27pt
$n_1$

\noindent\leftskip-11pt\vskip17pt
$n_2$

\leftskip5pt\vskip-20pt
$\theta$

\leftskip40pt\vskip2pt
$\theta$

\leftskip100pt\vskip-49pt

By [AGG, Lemma 4.2.2], the linear map $n_i:=\langle-,p\rangle u_i$, $v\mapsto\langle v,p\rangle u_i$,
corresponds, via the natural isomorphism in (2.1.1), to a normal vector to $P_i$ at $\pmb p$, $i=1,2$. Let
us show that $n_1$ points towards the half-space determined by $P_1$ containing $\pmb v_4$. Indeed, on one
hand, $\langle p,v_4\rangle=\langle v_1+v_2,v_4\rangle=g_{14}+g_{24}>0$. On~the other hand, by
Proposition 2.2.1,
$$\langle v_4,u_1\rangle\omega=\big\langle v_4,*(v_1\wedge v_2\wedge v_3)\big\rangle\omega=
v_4\wedge*\big(*v_1\wedge v_2\wedge v_3\big)=-v_1\wedge v_2\wedge v_3\wedge
v_4=-\sqrt{-\det G}\;\omega.$$
\leftskip0pt\noindent
This implies that $\langle v_4,u_1\rangle<0$ and, by Lemma 2.2.3, $n_1$ points towards the desired
direction. Similarly, one shows that $n_2$ points towards the half-space determined by $P_2$ {\it not\/}
containing $\pmb v_3$. Therefore (see the picture above),
$$\cos\angle(P_1,P_2)=\frac{\langle n_1,n_2\rangle}{\sqrt{\langle n_1,n_1\rangle}\sqrt{\langle
n_2,n_2\rangle}}=\frac{-\langle u_1,u_2\rangle}{\sqrt{-\langle u_1,u_1\rangle}\sqrt{-\langle
u_2,u_2\rangle}}.$$
It follows from Proposition 2.2.1 that
$$-\langle u_1,u_2\rangle=-\big\langle*(v_1\wedge v_2\wedge v_3),*(v_1\wedge v_2\wedge v_4)\big\rangle=\det
\left[\matrix0&&g_{12}&&g_{14}\\g_{21}&&0&&g_{24}\\g_{31}&&g_{32}&&g_{34}\endmatrix\right]=
g_{12}(g_{13}g_{24}+g_{14}g_{23}-g_{12}g_{34}),$$
$$-\langle u_1,u_1\rangle=-\big\langle*(v_1\wedge v_2\wedge v_3),*(v_1\wedge v_2\wedge v_3)\big\rangle=\det
\left[\matrix0&&g_{12}&&g_{13}\\g_{21}&&0&&g_{23}\\g_{31}&&g_{32}&&0\endmatrix\right]=2g_{12}g_{23}g_{31},$$
and
$$-\langle u_2,u_2\rangle=-\big\langle*(v_1\wedge v_2\wedge v_4),*(v_1\wedge v_2\wedge v_4)\big\rangle=\det
\left[\matrix0&&g_{12}&&g_{14}\\g_{21}&&0&&g_{24}\\g_{41}&&g_{42}&&0\endmatrix\right]=2g_{12}g_{24}g_{41}$$
which gives the desired formula.

Finally, if $\pmb v_4$ belongs to $P_1$, then the determinant of $G$ vanishes. This means that
$$g_{14}^2g_{23}^2+(g_{13}g_{24}-g_{12}g_{34})^2-2g_{14}g_{23}(g_{13}g_{24}+g_{12}g_{34})=0.$$
So,
$$\Big(\frac{g_{13}g_{24}+g_{14}g_{23}-g_{12}g_{34}}{2(g_{23}g_{31}g_{24}g_{41})^{\frac12}}
\Big)^2=\frac{g_{14}^2g_{23}^2+(g_{13}g_{24}-g_{12}g_{34})^2+2g_{14}g_{23}(g_{13}g_{24}-g_{12}g_{34})}
{4g_{23}g_{31}g_{24}g_{41}}=$$
$$=\frac{2g_{14}g_{23}(g_{13}g_{24}+g_{12}g_{34})+2g_{14}g_{23}(g_{13}g_{24}-g_{12}g_{34})}
{4g_{23}g_{31}g_{24}g_{41}}=1.$$
In other words,
$\displaystyle\frac{g_{13}g_{24}+g_{14}g_{23}-g_{12}g_{34}}{2(g_{23}g_{31}g_{24}g_{41})^{\frac12}}=\pm1$.
It is easy to verify (say, by considering particular cases in homogeneous coordinates), that the value $1$
occurs exactly when $\pmb v_3,\pmb v_4$ lie on a same half-plane determined by $\gamma$ and that the value
$-1$ occurs exactly when $\pmb v_3,\pmb v_4$ lie on distinct half-planes determined by $\gamma$.
\hfill$_\blacksquare$

\medskip

In the next section we will calculate the dihedral angle $\angle(P_1,P_2)$ in terms of some specific
representatives $v_1,v_2,v_3,v_4\in V$.

\bigskip

\centerline{\bf3.~Ideal tetrahedra}

\medskip

{\bf3.1.~Generalities.}~A {\it labelled ideal tetrahedron\/} is a $4$-tuple
$(\pmb v_1,\pmb v_2,\pmb v_3,\pmb v_4)$ of ideal points $\pmb v_i\in\S V$. Each $\pmb v_i$ is a {\it vertex\/} of
the labelled ideal tetrahedron. An {\it ideal tetrahedron\/} is a labelled ideal tetrahedron modulo permutations
of the vertices. A (labelled) ideal tetrahedron is {\it non-degenerate\/} when its vertices do not belong to a
same plane.

A {\it Gram matrix\/} of a labelled ideal tetrahedron $T=(\pmb v_1,\pmb v_2,\pmb v_3,\pmb v_4)$ is the Gram
matrix $G:=\big[\langle v_i,v_j\rangle\big]$ of representatives $v_1,v_2,v_3,v_4\in V$ of the vertices.
Clearly, two Gram matrices $G_1,G_2$ of $T$ are related by the expression $G_1=DG_2D$, where $D$ is a
diagonal matrix whose diagonal entries are non-null.

In this paper, we are mostly interested with Gram matrices that are doubly stochastic. (A {\it doubly
stochastic\/} matrix is a square matrix $A=[a_{ij}]$ with non-negative entries satisfying
$\sum_{i}{a_{ij}}=1$ for every $j$ and $\sum_{j}{a_{ij}}=1$ for every $i$.) A labelled ideal tetrahedron
admits at most one doubly stochastic Gram matrix (see Theorem 3.2.1). It turns out that the doubly
stochastic Gram matrix seems to be, in several circumstances, a natural choice among all the Gram matrices
of a given labelled ideal tetrahedron; these matrices are used in the formulation of Seidel's conjectures.

Let $\Cal T$ stand for the space of labelled ideal tetrahedra that admit a doubly stochastic Gram matrix
modulo isometries preserving the order of the vertices. As we will shortly see, $\Cal T$ is made up of the
classes of all labelled ideal tetrahedra with pairwise distinct vertices plus three particular classes of
degenerate labelled ideal tetrahedra. Only a few labelled ideal tetrahedra do {\it not\/} admit a doubly
stochastic Gram matrix; they are presented in Remark 3.2.3.

\bigskip

{\bf3.2.~Classifying triangles.}~The main result in this subsection is Theorem 3.2.1. It gives a particular
and explicit identification of the space $\Cal T$ with a plane equilateral triangle $\Delta$. In Proposition
3.2.5 we show that the symmetric $3$-group $S_3$ acts faithfully on $\Delta$ giving rise to the space of
(non-labelled) ideal tetrahedra modulo isometries.

\medskip

{\bf3.2.1.~Theorem.} {\it The space\/ $\Cal T$ can be identified with the equilateral triangle
$$\Delta:=\big\{(x,y,z)\in\Bbb R^3\mid x+y+z=1,x\leqslant y+z,y\leqslant z+x,z\leqslant x+y\big\}.$$
Explicitly, given a labelled ideal tetrahedron\/ $T=(\pmb v_1,\pmb v_2,\pmb v_3,\pmb v_4)$ admitting a
doubly stochastic Gram matrix, there exists a unique triple\/ $(r,s,t)\in\Delta$ such that the doubly
stochastic matrix
$$G:=\left[\matrix
0&r&s&t\\
r&0&t&s\\
s&t&0&r\\
t&s&r&0\endmatrix
\right]\leqno{(3.2.2)}$$
is a Gram matrix of\/ $T$ {\rm(}it is the only Gram matrix of\/ $T$ that is doubly stochastic\/{\rm)}. The
triple\/ $(r,s,t)\in\Delta$ depends only on the class of\/ $T$ modulo isometries preserving the order of the
vertices.

Conversely, given\/ $(r,s,t)\in\Delta$, there exists, up to isometry preserving the order of the vertices,
a~unique labelled ideal tetrahedron admitting the above doubly stochastic Gram matrix.

The interior of\/ $\Delta$ {\rm(}respectively, the boundary of\/ $\Delta${\rm)} corresponds to the
non-degenerate\/ {\rm(}respectively, the degenerate\/{\rm)} labelled ideal tetrahedra admitting a doubly
stochastic Gram matrix.}

\medskip

{\bf Proof.} Let $T=(\pmb v_1,\pmb v_2,\pmb v_3,\pmb v_4)$ be a labelled ideal tetrahedron with pairwise
distinct vertices and let $v_1,v_2,v_3,v_4$ be representatives of $\pmb v_1,\pmb v_2,\pmb v_3,\pmb v_4$. The
signature of the bilinear symmetric form $\langle-,-\rangle$ implies that $\langle v_i,v_j\rangle\neq0$ for
$i\neq j$. Rechoosing representatives, we assume that
$\langle v_1,v_2\rangle=\langle v_3,v_4\rangle>0$. Tha Gram matrix of $v_1,v_2,v_3$ has positive determinant
because $v_1,v_2,v_3$ span a space of signature $--+$. Hence,
$\langle v_1,v_3\rangle\langle v_2,v_3\rangle>0$. Simultaneously changing, if necessary, the signs of $v_1$
and $v_2$, we can assume that $\langle v_1,v_3\rangle,\langle v_2,v_3\rangle>0$. Considering the
determinants of the Gram matrices of $v_1,v_3,v_4$ and of $v_1,v_2,v_4$ we obtain
$\langle v_1,v_4\rangle,\langle v_2,v_4\rangle>0$. Now, scaling $v_1$ and $v_3$ by a same positive factor
allows us to consider $\langle v_1,v_3\rangle=\langle v_2,v_4\rangle$; an analogous reasoning involving
$v_1$ and $v_4$ leads to $\langle v_1,v_4\rangle=\langle v_2,v_3\rangle$. In~other words, $T$ admits
$G:=\left[\smallmatrix
0&r&s&t\\
r&0&t&s\\
s&t&0&r\\
t&s&r&0
\endsmallmatrix\right]$
as a Gram matrix. Dividing each representative $v_1,v_2,v_3,v_4$ by $\sqrt{r+s+t}$ we can assume that $G$
is doubly stochastic. Looking for a diagonal matrix $D$ whose diagonal entries are non-null and such that
$DGD$ is doubly stochastic, we obtain only two possibilities, namely $D=\pm I$. Hence, $G$ is the
unique Gram matrix of $T$ which is doubly stochastic.

Due to the signature of the form, $\det G=-(r+s+t)(-r+s+t)(r-s+t)(r+s-t)<0$ if $T$ is non-degenerate. This
means that $(-r+s+t),(r-s+t),(r+s-t)$ are positive numbers because, otherwise, one would obtain $r<0$ or
$s<0$ or $t<0$. Clearly, $\det G=0$ if $T$ is degenerate. Therefore, the obtained triple $(r,s,t)$ lies in
the interior of $\Delta$ if $T$ is non-degenerate and in the boundary of $\Delta$ if $T$ is degenerate. The
vertices of $\Delta$ are $\big(0,\frac12,\frac12\big)$, $\big(\frac12,0,\frac12\big)$,
$\big(\frac12,\frac12,0\big)$; so, $(r,s,t)$ cannot be a vertex of $\Delta$.

Regarding (degenerate) labelled ideal tetrahedra whose vertices are {\it not\/} pairwise distinct, it is not
difficult to see that only those of the forms $(\pmb v,\pmb v,\pmb v',\pmb v')$,
$(\pmb v,\pmb v',\pmb v,\pmb v')$, and $(\pmb v,\pmb v',\pmb v',\pmb v)$, $\pmb v\ne\pmb v'$, admit a doubly
stochastic Gram matrix. Such Gram matrices correspond respectively to the listed vertices of~$\Delta$.

Conversely, given $(r,s,t)\in\Delta$, there exists a labelled ideal tetrahedron admitting the doubly
stochastic Gram matrix $G$ by Sylvester's criterion.\footnote{Explicitly, if $(r,s,t)\in\Delta$ is not a
vertex of $\Delta$, then
$v_1:=b_1+b_2$, $v_2:=\frac{r}{2}b_1-\frac{r}{2}b_2$,
$v_3:=\left(\frac tr+\frac s2\right)b_1+\left(\frac tr-\frac s2\right)b_2+\sqrt{\frac{2st}r}b_3$,
$v_4:=\left(\frac sr+\frac t2\right)b_1+\left(\frac sr-\frac t2\right)b_2+
\frac{-r^2+s^2+t^2}{\sqrt{2rst}}b_3+\sqrt{\frac{-\det G}{2rst}}b_4$
(where $b_1,b_2,b_3,b_4$ is an orthonormal basis of $V$ satisfying
$\langle b_1,b_1\rangle=1$, $\langle b_2,b_2\rangle=\langle b_3,b_3\rangle=\langle b_4,b_4\rangle=-1$)
are the vertices of a labelled ideal tetrahedron having the required Gram matrix.}

Finally, it follows from [AGr1, Lemma 4.8.1] that two labelled ideal tetrahedra admit a same Gram matrix if
and only if they differ by an isometry preserving the order of the vertices.
$\hfill_\blacksquare$

\medskip

{\bf3.2.3.~Remark.}~The labelled ideal tetrahedra that do not admit a doubly stochastic Gram matrix are the
following: exactly $2$ vertices coincide (as in $(\pmb v,\pmb v,\pmb v',\pmb v'')$, $v''\ne v\ne v'\ne v''$,
for example) or at least $3$ vertices coincide (as in $(\pmb v,\pmb v,\pmb v,\pmb v')$, for example). In
what follows, we only deal with labelled ideal tetrahedra that admit a doubly stochastic Gram matrix. So,
whenever we refer to labelled ideal tetrahedra, we are actually referring to those that admit a doubly
stochastic Gram matrix. The same goes for (non-labelled) ideal tetrahedra, as we only consider those arising
from labelled ideal tetrahedra that admit a doubly stochastic Gram matrix.
\hfill$_\blacksquare$

\medskip

In view of the identification $\Delta\simeq\Cal T$ we will also refer to a point in $\Delta$ as a `labelled
ideal tetrahedron.'

\medskip

{\bf3.2.4.~Remark.}~It is curious to note that $\Delta$ is nothing but the space of Euclidean triangles
with ordered vertices and fixed perimeter ($=1$) modulo isometries of the plane that preserve the order of
the vertices.

Let $(r,s,t)\in\Delta$ be a point that is not a vertex of $\Delta$. Then $(r,s,t)$ corresponds to an
Euclidean triangle with ordered vertices whose vertices are pairwise distinct. Let
$(\theta_1,\theta_2,\theta_3)$ be the interior angles of this Euclidean triangle\footnote{Here, we exclude
the degenerate Euclidean triangles corresponding to the vertices of $\Delta$ because, in this case, two
vertices of the Euclidean triangle coincide and the internal angles are not well-defined.} (the angle
$\theta_1$ is opposite to a side of length $r$, the angle $\theta_2$ is opposite to a side of length $s$,
and the angle $\theta_3$ is opposite to a side of length $t$). Hence,
$$\cos\theta_1:=\frac{-r^2+s^2+t^2}{2st},\quad
\cos\theta_2:=\frac{r^2-s^2+t^2}{2rt},\quad
\cos\theta_3:=\frac{r^2+s^2-t^2}{2rs}.$$
Applying Proposition 2.2.4 to the doubly stochastic Gram matrix (3.2.2) of the labelled ideal tetrahedron
$T=(r,s,t)\in\Delta$, it is straightforward to see that the angles $\theta_1,\theta_2,\theta_3$ are
precisely the dihedral angles of~$T$. More specifically, the dihedral angles at the edge joining the
vertices $v_1,v_2$ and at the edge joining $v_3,v_4$ are both $\theta_1$; the dihedral angles at the edge
joining the vertices $v_1,v_3$ and at the edge joining $v_2,v_4$ are both~$\theta_2$; the dihedral angles at
the edge joining the vertices $v_1,v_4$ and at the edge joining $v_2,v_3$ are both $\theta_3$.

In particular, opposite dihedral angles of an ideal tetrahedron are equal and dihedral angles incident to a
same vertex sum $\pi$~---~a couple of well-known facts (see, for instance, [Mil] or [Sei]).
\hfill$_\blacksquare$

\medskip

In order to obtain the space of (non-labelled) ideal tetrahedra modulo isometries, let us study the action of the
symmetric $4$-group $S_4$ on $\Delta$ by permutations of vertices of labelled ideal tetrahedra.

\medskip

{\bf3.2.5.~Proposition.}~{The kernel of the natural action of the symmetric $4$-group\/ $S_4$ on\/ $\Delta$
is isomorphic to the Klein four group\/ $H$. The symmetric $3$-group $S_3=S_4/H$ acts on\/ $\Delta$ by
permutations of coordinates or, equivalently, by reflections on the altitudes of\/ $\Delta$.}

\medskip

{\bf Proof.}~Clearly, $S_4$ acts on $\Delta$ and every permutation of coordinates in $\Delta$ is induced by
the action of an element of $S_4$. Let $p\in H$ be a permutation in the Klein four group
$$H:=\big\{I,(12)(34),(13)(24),(14)(23)\big\}\leqslant S_4.$$
A direct calculation shows that the labelled
ideal tetrahedra $(\pmb v_1,\pmb v_2,\pmb v_3,\pmb v_4)$ and
$(\pmb v_{p(1)},\pmb v_{p(2)},\pmb v_{p(3)},\pmb v_{p(4)})$ have the same doubly stochastic Gram matrix
(3.2.2). Hence, $H$ is contained in the kernel of the $S_4$-action. Taking a point $(r,s,t)\in\Delta$
with pairwise distinct $r,s,t$, there are six points in $\Delta$ corresponding to the permutations of
$r,s,t$ (they are the reflections of $(r,s,t)$ on the altitudes $x=y$, $y=z$, and $z=x$ of $\Delta$) and,
therefore, the action of $S_3=S_4/H$ on $\Delta$ is faithful.
\hfill$_\blacksquare$

\medskip

By Remark 3.2.4, points in the altitudes of $\Delta$ (except the vertices of $\Delta$) correspond to the
`isosceles' tetrahedra (i.e., two dihedral angles are equal); the tetrahedron $r=s=t=\frac13$ is the
`equilateral' or `regular' one (i.e., all of its dihedral angles are equal). These tetrahedra have fewer
copies by the $S_3$ action on $\Delta$ because, in comparison to the generic case, there are more
permutations of their vertices that can be achieved by means of isometries of the ambient space (after all,
they are more regular). In particular, in the equilateral case, every permutation of vertices can be
accomplished through an isometry.

\smallskip

The altitudes of $\Delta$ divide this triangle into six congruent triangles $\Delta_i$, $i=1,\dots,6$; we
put
$$\Delta_1:=\big\{(x,y,z)\in\Delta\mid x\leqslant y\leqslant z\big\}.\leqno{(3.2.6)}$$

{\bf3.2.7.~Corollary.} {\it Each\/ $\Delta_i$ is a copy of the space of\/ {\rm(}non-labelled\/{\rm)} ideal
tetrahedra modulo isometries.}

\medskip

In order to prove Seidel's fourth conjecture it will be easier to take, in place of $\Delta$, an equilateral
triangle centred at the origin of a plane:

\medskip

{\bf3.2.8.~Remark.} We identify the affine plane $\big\{(x,y,z)\mid x+y+z=1\big\}\subset\Bbb R^3$ with
$\Bbb R^2$ by taking $(\frac13,\frac13,\frac13)$ as the origin and
$\left(0,-\frac12,\frac12\right),\frac{\sqrt3}3\left(-1,\frac12,\frac12\right)$ as a basis. In these
coordinates, $\Delta$ becomes the equilateral triangle
$$\widetilde{\Delta}:=\left\{(c,d)\in\Bbb R^2\mid d\geqslant-\frac{\sqrt3}6,d\leqslant-\sqrt3c+
\frac{\sqrt3}3,d\leqslant\sqrt3c+\frac{\sqrt3}3\right\}\subset\Bbb R^2$$
centred at the origin. The point $(c,d)\in\widetilde\Delta$ corresponds to the point
$\left(\frac{1-\sqrt{3}d}3,\frac{2-3c+\sqrt3d}6,\frac{2+3c+\sqrt3d}6\right)\in\Delta$. Each $\Delta_i$
will be denoted by $\widetilde\Delta_i$ in the new coordinates; in particular, the triangle $\Delta_1$ in
(3.2.6) becomes
$$\widetilde{\Delta}_1:=\left\{(c,d)\in\Bbb R^2\mid c\geqslant0,d\geqslant\frac{\sqrt3}3\cdot c,d\leqslant-
\sqrt3c+\frac{\sqrt3}3\right\}.\eqno{\hfill_\blacksquare}$$

\centerline{\bf4.~Proof of Seidel's first conjecture}

\medskip

We now show that determinant and permanent of doubly stochastic Gram matrices of labelled ideal tetrahedra
can be taken as coordinates of the space of ideal tetrahedra modulo isometries. In other words, these
algebraic functions uniquely determine an ideal tetrahedron modulo isometries. This is a stronger version of
Seidel's first conjecture which claims that volume is a function of determinant and permanent.

\smallskip

We remind the reader that the permanent of an $m\times m$ matrix $A=[a_{ij}]$ is defined by the expression
$$\per A:=\sum_{\sigma\in S_m}a_{1\sigma(1)}a_{2\sigma(2)}\dots a_{m\sigma(m)},$$
where $S_m$ stands for the symmetric $m$-group.

\medskip

{\bf4.1.~Theorem [{\rm Seidel's Speculation 1}].} {\it Let\/ $\alpha,\beta\in\Bbb R$ be respectively the
determinant and permanent of the doubly stochastic Gram matrix of a labelled ideal tetrahedron\/
$T\in\Delta$. Then the only tetrahedra in\/~$\Delta$ whose doubly stochastic Gram matrices have determinant
and permanent respectively equal to $\alpha,\beta$ are those that differ from\/ $T$ by the\/ $S_3$-action in
Proposition\/ {\rm3.2.5.} In other words, the pair\/ $(\alpha,\beta)$ determines a unique\/
{\rm(}non-labelled\/{\rm)} ideal tetrahedron.}

\medskip

{\bf Proof.} Let
$G:=\left[\smallmatrix
0&r&s&t\\
r&0&t&s\\
s&t&0&r\\
t&s&r&0\endsmallmatrix
\right]$
be the doubly stochastic Gram matrix of $T$. Then $r,s,t$ satisfy the equations
$$\cases x+y+z=1\\
-(x+y+z)(-x+y+z)(x-y+z)(x+y-z)=\alpha\\
(x^2+y^2+z^2)^2=\beta\\
\endcases\leqno{(4.2)}$$
which are equivalent to
$$\cases x+y+z=1\\
xy+yz+zx=\frac{1-\sqrt\beta}2\\
xyz=\frac{\alpha-2\sqrt\beta+1}8\endcases.\leqno{(4.3)}$$
So, $r,s,t$ are the roots of the polynomial
$$p(w):=(w-r)(w-s)(w-t)=w^3-(r+s+t)w+(rs+st+tr)w-rst=
w^3-w^2+\frac{1-\sqrt\beta}2w-\frac{\alpha-2\sqrt\beta+1}8.$$
If a triple $(r',s',t')\in\Delta$ satisfies $(4.3)$, then $r',s',t'$ are the roots of the polynomial
$(w-r')(w-s')(w-t')=p(w)$, that is, $(r',s',t')$ is a permutation of $(r,s,t)$.
$\hfill_\blacksquare$

\medskip

Given a tetrahedron $T\in\Delta$, let $G_T$ denote the doubly stochastic Gram matrix of $T$. We define
$$S:=\Big\{(\alpha,\omega)\mid\alpha=\det G_T\text{\rm\;\, and }\omega=\sqrt{\per G_T}
\text{\rm\;\,for }T\in\Delta\Big\}\subset\Bbb R^2.\leqno{(4.4)}$$

Let $\Delta_i$ and $\widetilde\Delta_i$ be as defined in the previous section (see Remark 3.2.8). We have the
following corollary to Theorem 4.1:

\medskip

{\bf4.5.~Corollary.} {\it The function\/ $\Delta_i\to S$, $T\mapsto\big(\det G_T,\sqrt{\per G_T}\big)$, is a
homeomorphism for each\/ $i=1,2,\dots,6$. {\rm(}Of course, the same is true for $\widetilde\Delta_i$ in
place of $\Delta_i$.{\rm)}}

\medskip

{\bf Proof.} The map in question is clearly continuous; it is bijective by Theorem 4.1. Since $\Delta_i$ is
compact and $S$ is Hausdorff, the map is a homeomorphism.
$\hfill_\blacksquare$

\medskip

$\epsfxsize=3cm\vcenter{\hbox{\epsfbox{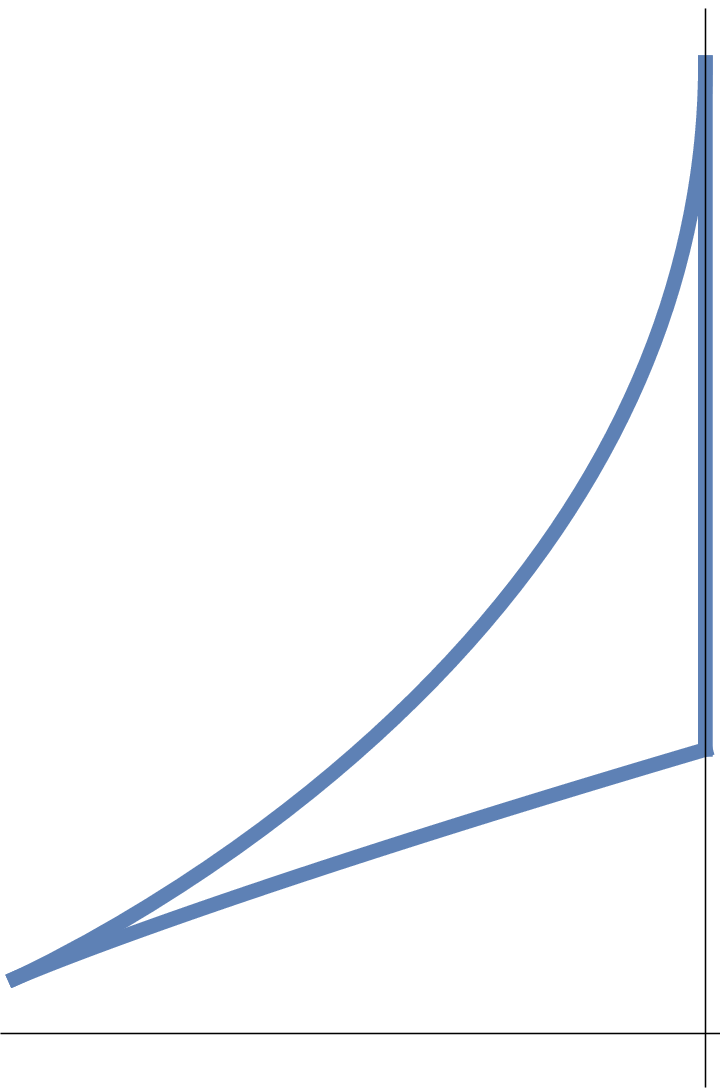}}}$

\vskip-40pt\noindent
$(-\frac1{27},\frac13)$

\leftskip90pt\vskip-20pt
$(0,\frac38)$

\vskip-90pt
$(0,\frac12)$

\vskip55pt\leftskip65pt
$S$

\vskip40pt\leftskip60pt
$\det$

\vskip-90pt\leftskip88pt
$c_3$

\leftskip48pt
$c_1$

\leftskip45pt\vskip30pt
$c_2$

\leftskip85pt\vskip-12pt
$\sqrt{\per}$

\vskip-115pt\leftskip140pt\noindent
The space $S$, depicted on the left, can be seen as a reparameterization of the space $\Delta_i$ of
(non-labelled) ideal tetrahedra modulo isometries; in this reparameterization, the determinant and (the
square root of the) permanent of doubly stochastic Gram matrices of labelled ideal tetrahedra are the
coordinates. The space $S$ will be used in the proof of Seidel's fourth conjecture (see Sections 5 and 6)
and it is explicitly described in the proof of Theorem 6.10. The curves $c_1$ and $c_2$ correspond to the
isosceles ideal tetrahedra; the curve $c_3$ lists the degenerate ideal tetrahedra. The vertices
$(-\frac1{27},\frac13)$, $(0,\frac12)$, and $(0,\frac38)$ correspond respectively to the equilateral (or
regular) ideal tetrahedron, to the degenerate ideal tetrahedron with two pairs of coinciding vertices, and
to the `regular' degenerate ideal tetrahedron.\footnote{The latter can be described as follows. Take four
pairwise distinct points $\pmb v_1,\pmb v_2,\pmb v_3,\pmb v_4$ on the ideal boundary of a same plane.
Assume that these points are oriented in the counterclockwise sense. Then the geodesic joining
$\pmb v_1,\pmb v_3$ and the geodesic joining $\pmb v_2,\pmb v_4$ intersect orthogonally.}

\leftskip0pt

\bigskip

\centerline{\bf5.~A volume formula}

\medskip

\noindent
The volume of an ideal tetrahedron is the hyperbolic volume of the convex hull of its vertices. In this
section we provide an explicit formula for the volume of an ideal tetrahedron in terms of the coordinates of
the space $S$ of ideal tetrahedra modulo isometries (see Corollaries 3.2.7 and 4.5). In other words,
Theorem~5.1 describes the volume of an ideal tetrahedron as a function of the determinant and the (square
root of the) permanent of doubly stochastic Gram matrices of labelled ideal tetrahedra.

\smallskip

Let $T$ be an ideal tetrahedron with dihedral angles $\theta_1,\theta_2,\theta_3$ (see Remark 3.2.4). Milnor's
beautiful volume formula [Mil] states that the volume of $T$ is given by
$$\vol(T)=\lob(\theta_1)+\lob(\theta_2)+\lob(\theta_3),$$
where $\lob:\Bbb R\to\Bbb R$ is the Lobachevsky function defined in (1.1). This is the volume formula we are
going to use in the next theorem.

\medskip

{\bf5.1.~Theorem.} {\it  The volume function\/ $\vol:S\to\Bbb R$ is
given by\/
$$\vol(\alpha,\omega)=\cases\lob(\theta_1)+\lob(\theta_2)+\lob(\theta_3),&(\alpha,\omega)\ne
\left(0,\frac12\right)\\0,&(\alpha,\omega)=\left(0,\frac12\right)\endcases,$$
where\/
$$\theta_1:=\arccos\frac{-r^2+s^2+t^2}{2st},\quad
\theta_2:=\arccos\frac{r^2-s^2+t^2}{2rt},\quad
\theta_3:=\arccos\frac{r^2+s^2-t^2}{2rs},$$
$$(r,s,t):=\left(\frac{1-\sqrt{3}d}{3},\frac{2-3c+\sqrt{3}d}{6},\frac{2+3c+\sqrt{3}d}{6}\right),$$
and
$$(c,d):=\cases
(0,0),&(\alpha,\omega)=\left(-\frac{1}{27},\frac{1}{3}\right)\\
\sqrt{2\omega-\frac23}\cdot(\sin\kappa,\cos\kappa),&(\alpha,\omega)\neq
\left(-\frac1{27},\frac13\right)\endcases,\quad
\kappa:=\displaystyle\frac13\arccos\frac{-27\alpha+18\omega-7}
{4\sqrt2\left(3\omega-1\right)^{\frac32}}.$$}

\medskip

{\bf Proof.} Let $(\alpha,\omega)\in S$. Writing the equations (4.2) in terms of the coordinates
$(c,d)\in\widetilde\Delta$ (see Remark 3.2.8) one obtains
$$\cases
\displaystyle\frac{\left(2\sqrt3d+1\right)\left(9c^2-\left(\sqrt3d-1\right)^2\right)}{27}=\alpha\\
\displaystyle\frac{3c^2+3d^2+2}6=\omega\\
\endcases.\leqno{(5.2)}$$
By Theorem 4.1, these equations have a unique solution in $\widetilde\Delta_1$. If
$(\alpha,\omega)=(-\frac1{27},\frac13)$, it is easy to see that this solution is
$(0,0)\in\widetilde\Delta_1$. Let us assume that $(\alpha,\omega)\ne(-\frac1{27},\frac13)$.

It follows from (5.2) that $2\sqrt3d$ satisfies the cubic equation $x^3+a_0x+b_0=0$ with
$$a_0:=6(1-3\omega),\quad b_0:=27\alpha-18\omega+7.\leqno{(5.3)}$$
A straightforward calculation shows that the cubic equation $x^3+ax+b=0$ has the three roots
$$x_k:=2\sqrt{\frac{-a}3}\cos\left(\frac13\arccos\frac{-3\sqrt3b}{2\sqrt{-a^3}}+\frac{2k\pi}{3}\right),
\quad k=0,1,2,$$
when $a,b\in\Bbb R$ are such that $a\neq 0$ and $\displaystyle\frac{b^2}4+\frac{a^3}{27}\leqslant0$. Let us
show that $a_0,b_0$ satisfy the previous inequalities. Indeed,
$$a_0=6\left(1-3\frac{3c^2+3d^2+2}6\right)=-9\left(c^2+d^2\right)\ne0$$
because $(c,d)=(0,0)$ leads to $(\alpha,\omega)=(-\frac1{27},\frac13)$. Moreover,
$$b_0=27\frac{\left(2\sqrt3d+1\right)\left(9c^2-\left(\sqrt3d-1\right)^2\right)}{27}-18
\frac{3c^2+3d^2+2}{6}+7=6\sqrt3d\left(3c^2-d^2\right)$$
implies
$$\frac{b_0^2}{4}+\frac{a_0^3}{27}=-27c^2\left(3d^2-c^2\right)^2\leqslant0.\leqno{(5.4)}$$
Therefore, taking $d_k:=\frac1{2\sqrt3}\cdot x_k$ and defining $c_k:=\sqrt{2\omega-d_k^2-\frac23}$,
$k=0,1,2$, provides all possible solutions $(\pm c_k,d_k)\in\widetilde\Delta$ of equations (5.2). The
solution $(c_0,d_0)$ lies in $\widetilde\Delta_1$ since $0\leqslant\kappa\leqslant\frac\pi3$ implies
$c_0\geqslant0$ and $d_0\geqslant\frac{\sqrt3}3\cdot c_0$.

It remains to apply Remarks 3.2.8, 3.2.4, and Milnor's volume formula.
\hfill$_\blacksquare$

\medskip

It is well-known that the regular ideal tetrahedron (the one whose dihedral angles are all equal to~$\pi/3$)
has maximal volume among all tetrahedra in hyperbolic $3$-space [Mil]. This is the tetrahedron corresponding
to the point $(\alpha,\omega)=(-1/27,1/3)$ (see also the explicit description of $S$ in the proof of Theorem
6.10).

\bigskip

\centerline{\bf6.~Proof of Seidel's fourth conjecture}

\medskip

This section contains the main result of the paper. In Theorem 6.10 it is shown that the volume function
$\vol:S\to\Bbb R$ from the space $S$ of ideal tetrahedra modulo isometries (see Corollaries 3.2.7 and~4.5)
is monotonic both in the determinant and in the permanent of doubly stochastic Gram matrices of labelled
ideal tetrahedra. Theorem 6.10 follows from a long calculation (presented below) which had to be approached
quite carefully. Indeed, in many places, we had to make some appropriate choices of coordinates in order to
keep the expressions treatable.

It is curious to note that, at the end of the day, proving that the volume is a monotonic function of the
determinant and of the permanent of doubly stochastic Gram matrices of labelled ideal polyhedra amounts
to determining the sign of a {\it single\/} expression, the term $k'$ defined in (6.9). Even more intriguing
is the fact that studying the sign of $k'$ is, in some sense, the most delicate part of the proof of the
fourth conjecture.

\medskip

Let $O:=\big\{(r,s,t)\in\Bbb R^3\mid r<s<t;\,r<s+t;\,s<r+t;\,t<r+s\big\}$ be an open set in $\Bbb R^3$.
Using the notation from Theorem 5.1, we define the functions
$$h:\overset{\;\circ}\to S\to\widetilde\Delta_1,\quad
(\alpha,\omega)\mapsto\sqrt{2\omega-\frac23}\cdot(\sin\kappa,\cos\kappa)\text{\rm,\; where }
\kappa:=\displaystyle\frac13\arccos\frac{-27\alpha+18\omega-7}
{4\sqrt2\left(3\omega-1\right)^{\frac32}},$$
$$g:\overset{\!\!\!\circ}\to{\widetilde\Delta_1}\to\Delta_1,\quad(c,d)\mapsto
\left(\frac{1-\sqrt{3}d}{3},\frac{2-3c+\sqrt{3}d}{6},\frac{2+3c+\sqrt{3}d}{6}\right),$$
and
$$f:O\to\Bbb R,\quad(r,s,t)\mapsto
\lob(\theta_1)+\lob(\theta_2)+\lob(\theta_3),$$
where
$$\theta_1:=\arccos\frac{-r^2+s^2+t^2}{2st},\quad
\theta_2:=\arccos\frac{r^2-s^2+t^2}{2rt},\quad
\theta_3:=\arccos\frac{r^2+s^2-t^2}{2rs}.$$
Theorem 5.1 implies that the restriction $\vol:\overset{\;\circ}\to S\to\Bbb R$ is given by
$\vol=f\circ g\circ h$.

Let $h_1,h_2$ and $g_1,g_2,g_3$ stand respectively for the coordinate functions of $h$ and $g$. In what
follows, we will calculate the derivatives
$$\frac{\partial\vol}{\partial\alpha}=\left(\frac{\partial f}{\partial r}\frac{\partial g_1}{\partial
c}+\frac{\partial f}{\partial s}\frac{\partial g_2}{\partial c}+\frac{\partial f}{\partial t}\frac{\partial
g_3}{\partial c}\right)\frac{\partial h_1}{\partial\alpha}+\left(\frac{\partial f}{\partial r}\frac{\partial
g_1}{\partial d}+\frac{\partial f}{\partial s}\frac{\partial g_2}{\partial d}+\frac{\partial f}{\partial
t}\frac{\partial g_3}{\partial d}\right)\frac{\partial h_2}{\partial\alpha}=\leqno{(6.1)}$$
$$=\frac12\left(-\frac{\partial f}{\partial s}+\frac{\partial f}{\partial t}\right)\frac{\partial h_1}
{\partial\alpha}+\frac{\sqrt3}6\left(-2\frac{\partial f}{\partial r}+\frac{\partial f}{\partial
s}+\frac{\partial f}{\partial t}\right)\frac{\partial h_2}{\partial\alpha}$$
and
$$\frac{\partial\vol}{\partial\omega}=\frac12\left(-\frac{\partial f}{\partial s}+\frac{\partial f}{\partial
t}\right)\frac{\partial h_1}{\partial\omega}+\frac{\sqrt3}6\left(-2\frac{\partial f}{\partial
r}+\frac{\partial f}{\partial s}+\frac{\partial f}{\partial t}\right)\frac{\partial
h_2}{\partial\omega}\leqno{(6.2)}$$
at each point $(\alpha,\omega)\in\overset\circ\to S$.

\medskip

{\bf6.3.~Lemma.}~{\it At each point\/ $(\alpha,\omega)\in\overset{\;\circ}\to S$ we have
$$\frac{\partial h_1}{\partial\alpha}=
\frac{\sqrt3h_2}{2h_1\left(3h_2^2-h_1^2\right)},\quad
\frac{\partial h_2}{\partial\alpha}=-\frac{\sqrt3}{2\left(3h_2^2-h_1^2\right)}$$
$$\frac{\partial h_1}{\partial\omega}=
\frac{3h_2^2-3h_1^2-\sqrt3h_2}{3h_1\left(3h_2^2-h_1^2\right)},\quad
\frac{\partial h_2}{\partial\omega}=
\frac{6h_2+\sqrt3}{3\left(3h_2^2-h_1^2\right)}.$$}

\medskip

{\bf Proof.} Note that
$\displaystyle\frac{\partial h_1}{\partial\alpha}=\frac{3\sqrt{6(3\omega-1)}}l\cos\kappa=
\frac9l\cdot h_2$,
where
$$l:=\sqrt{\left(4\sqrt2\left(3\omega-1\right)^{\frac32}\right)^2-\left(-27\alpha+18\omega-7\right)^2}.$$
Let $a_0,b_0$ be as in (5.3). Then
$$l=\sqrt{32(3\omega-1)^3-(-27\alpha+18\omega-7)^2}=\sqrt{32\left(-\frac{a_0}6\right)^3-(-b_0)^2}=
2\sqrt{-\left(\frac{b_0^2}4+\frac{a_0^3}{27}\right)}$$
because $\big(h_1(\alpha,\omega),h_2(\alpha,\omega)\big)$ satisfies equations (5.2). So,
(5.4) implies that $l=2\sqrt{27h_1^2\left(3h_2^2-h_1^2\right)^2}$. It~follows from
$\big(h_1(\alpha,\omega),h_2(\alpha,\omega)\big)\in\overset{\!\!\!\circ}\to{\widetilde\Delta_1}$
that $h_1>0$ and $3h_2^2>h_1^2$; hence,
$$l=6\sqrt3h_1(3h_2^2-h_1^2).\leqno{(6.4)}$$
The expression for $\displaystyle\frac{\partial h_2}{\partial\alpha}$ is obtained in a similar fashion.

Concerning the derivatives with respect to $\omega$, we have
$$\frac{\partial h_1}{\partial\omega}=\frac{\sqrt3}{\sqrt2\sqrt{3\omega-1}}\sin\kappa-
\frac{3\sqrt3(9\alpha-2\omega+1)}{l\sqrt2\sqrt{3\omega-1}}\cos\kappa=\frac3{2(3\omega-1)}h_1
-\frac{9(9\alpha-2\omega+1)}{2(3\omega-1)l}h_2.$$
But $2(3\omega-1)=3(h_1^2+h_2^2)$ and
$$9\alpha-2\omega+1=9\frac{\left(2\sqrt3h_2+1\right)\left(9h_1^2-\left(\sqrt3h_2-1\right)^2\right)}{27}-
2\frac{3h_1^2+3h_2^2+2}{6}+1=$$
$$=2\left(3\sqrt3h_1^2h_2+h_1^2-\sqrt3h_2^3+h_2^2\right)$$
because $\big(h_1(\alpha,\omega),h_2(\alpha,\omega)\big)$ satisfies equations (5.2). It remains to use
the expression for $l$ in (6.4). The derivative $\displaystyle\frac{\partial h_2}{\partial\omega}$ is
obtained analogously.
\hfill$_\blacksquare$

\medskip

{\bf6.5.~Lemma.}~{\it At each point\/
$(r,s,t)\in\overset\circ\to\Delta_1=\big\{(r,s,t)\in O\mid r+s+t=1\big\}$ we have
$$\frac{\partial f}{\partial r}=\frac2{\sqrt{-\alpha}}(-r\log r+\cos\theta_3s\log s+\cos\theta_2t\log t),$$
$$\frac{\partial f}{\partial s}=\frac2{\sqrt{-\alpha}}(\cos\theta_3r\log r-s\log s+\cos\theta_1t\log t),$$
$$\frac{\partial f}{\partial t}=\frac2{\sqrt{-\alpha}}(\cos\theta_2r\log r+\cos\theta_1s\log s-t\log t).$$}

\medskip

{\bf Proof.} A straightforward calculation shows that, at each point $(r,s,t)\in O$,
$$\frac{\partial}{\partial r}\lob(\theta_1)=-\frac{2r}k\log\left(\frac k{st}\right),\quad
\frac{\partial}{\partial r}\lob(\theta_2)=\frac{r^2+s^2-t^2}{rk}\log\left(\frac k{rt}\right),\quad
\frac{\partial}{\partial r}\lob(\theta_3)=\frac{r^2-s^2+t^2}{rk}\log\left(\frac k{rs}\right),$$
where $k:=\sqrt{(r+s+t)(-r+s+t)(r-s+t)(r+s-t)}$. Adding these equations, one obtains
$$\frac{\partial f}{\partial r}=\frac1k\left(-2r\log r+\frac{r^2+s^2-t^2}r\log
s+\frac{r^2-s^2+t^2}r\log t\right)=
\frac2k(-r\log r+\cos\theta_3s\log s+\cos\theta_2t\log t).$$
The second equation in (4.2) says that $k=\sqrt{-\alpha}$ for $(r,s,t)\in\overset\circ\to\Delta_1$.
The remaining derivatives follow by symmetry.
\hfill$_\blacksquare$

\medskip

For notational simplicity, in what follows, we will write $c=h_1(\alpha,\omega)$, $d=h_2(\alpha,\omega)$,
$r=g_1\big(h(\alpha,\omega)\big)$, $s=g_2\big(h(\alpha,\omega)\big)$, and $t=g_3\big(h(\alpha,\omega)\big)$.

{\bf6.6.~Proposition.}~{\it At each point\/ $(\alpha,\omega)\in\overset{\;\circ}\to S$ we have
$$\frac{\partial\vol}{\partial\alpha}=\frac{\sqrt3}{648c\left(3d^2-c^2\right)rst\sqrt{-\alpha}}
\left(M\log\frac{st}{r^2}+N\log\frac{t}{s}\right),$$
where
$$M:=c(3d-\sqrt3)(21d^2+4\sqrt3d+9c^2-2),$$
$$N:=\sqrt3\big(-27c^4+9c^2(3d^2+\sqrt3d+1)+d(18d^3+9\sqrt3d^2-9d-2\sqrt3)\big).$$
Moreover,
$$\frac{\partial\vol}{\partial\omega}=\frac{\sqrt{-\alpha}}
{12c\left(3d^2-c^2\right)rst}\left(P\log\frac{st}{r^2}+Q\log\frac ts\right),$$
where
$$P:=2c(\sqrt3d-1)$$
$$Q:=-3c^2+d(3d+2\sqrt3).$$}

\medskip

{\bf Proof.}~By Equations (6.1) and (6.2),
$$\frac{\partial\vol}{\partial\alpha}=\frac12\left(-\frac{\partial f}{\partial s}+\frac{\partial f}{\partial
t}\right)\frac{\partial h_1}{\partial\alpha}+\frac{\sqrt3}6\left(-2\frac{\partial f}{\partial
r}+\frac{\partial f}{\partial s}+\frac{\partial f}{\partial t}\right)\frac{\partial h_2}{\partial\alpha}$$
and
$$\frac{\partial\vol}{\partial\omega}=\frac12\left(-\frac{\partial f}{\partial s}+\frac{\partial f}{\partial
t}\right)\frac{\partial h_1}{\partial\omega}+\frac{\sqrt3}6\left(-2\frac{\partial f}{\partial
r}+\frac{\partial f}{\partial s}+\frac{\partial f}{\partial t}\right)\frac{\partial
h_2}{\partial\omega}.$$
On one hand, it follows from Lemma 6.5 that
$$\frac12\left(-\frac{\partial f}{\partial s}+\frac{\partial f}{\partial t}\right)=\frac1{\sqrt{-\alpha}}
\Bigg(\bigg(-\frac{1-2t}{2s}+\frac{1-2s}{2t}\bigg)\log r+\frac{1-2r}{2t}\log s-\frac{1-2r}{2s}\log
t\Bigg)$$
since
$$\cos\theta_1=\frac{1-2r}{2st}-1,\quad\cos\theta_2=\frac{1-2s}{2rt}-1,\quad
\cos\theta_3=\frac{1-2t}{2rs}-1$$
due to $r+s+t=1$.
So,
$$\frac12\left(-\frac{\partial f}{\partial s}+\frac{\partial f}{\partial t}\right)=
\frac{A\log r+B\log s+C\log t}{2rst\sqrt{-\alpha}}$$
where $A:=r\big(s(1-2s)-t(1-2t)\big)$, $B:=rs(1-2r)$, $C:=-rt(1-2r)$.
The equalities $r=\frac{1-\sqrt{3}d}{3}$, $s=\frac{2-3c+\sqrt3d}6$, and $t=\frac{2+3c+\sqrt3d}6$ imply
$$A=\frac c9\left(-6d^2+\sqrt3d+1\right),$$
$$B=\frac1{54}\left(18cd^2-3\sqrt3cd-3c-6\sqrt3d^3-9d^2+3\sqrt3d+2\right),$$
$$C=\frac1{54}\left(18cd^2-3\sqrt3cd-3c+6\sqrt3d^3+9d^2-3\sqrt3d-2\right).$$
On the other hand,
$$\frac{\sqrt3}6\left(-2\frac{\partial f}{\partial r}+\frac{\partial f}{\partial s}+\frac{\partial
f}{\partial t}\right)=\frac{\sqrt3}{3\sqrt{-\alpha}}\Bigg(\bigg(\frac{1-2t}{2s}+\frac{1-2s}{2t}\bigg)\log
r+\bigg(-\frac{1-2t}r+\frac{1-2r}{2t}\bigg)\log s+$$
$$+\bigg(-\frac{1-2s}r+\frac{1-2r}{2s}\bigg)\log t\Bigg)=
\frac{D\log r+E\log s+F\log t}{6rst\sqrt{-\alpha}},$$
where $D:=\sqrt3r\big(t(1-2t)+s(1-2s)\big)$, $E:=\sqrt3s\big(r(1-2r)-2t(1-2t)\big)$,
$F:=\sqrt3t\big(r(1-2r)-2s(1-2s)\big)$.
From above expressions for $r,s,t$ one obtains
$$D=\frac1{27}\left(27c^2d-9\sqrt3c^2+9d^3-9d+2\sqrt3\right),$$
$$E=\frac1{54}\left(-27\sqrt3c^3-27c^2d+9\sqrt3c^2+27\sqrt3cd^2+27cd+9\sqrt3c-9d^3+9d-2\sqrt3\right),$$
$$F=\frac1{54}\left(27\sqrt3c^3-27c^2d+9\sqrt3c^2-27\sqrt3cd^2-27cd-9\sqrt3c-9d^3+9d-2\sqrt3\right).$$
Therefore, by Lemma 6.3,
$$\frac{\partial\vol}{\partial\alpha}=
\frac{(3dA-cD)\log r+(3dB-cE)\log s+(3dC-cF)\log t}{4\sqrt3c(3d^2-c^2)rst\sqrt{-\alpha}}=$$
$$=\frac{\sqrt3}{12c(3d^2-c^2)rst\sqrt{-\alpha}}\left(-\frac{M}{27}\log r+\frac{M-N}{54}\log
s+\frac{M+N}{54}\log t\right)=$$
$$=\frac{\sqrt3}{648c\left(3d^2-c^2\right)rst\sqrt{-\alpha}}
\left(M\log\frac{st}{r^2}+N\log\frac{t}{s}\right).$$
Finally,
$$\frac{\partial\vol}{\partial\omega}=\frac{
\big(3(3d^2-3c^2-\sqrt3d)A+c(6d+\sqrt3)D\big)\log
r+\big(3(3d^2-3c^2-\sqrt3d)B+c(6d+\sqrt3)E\big)\log
s}{18c(3d^2-c^2)rst\sqrt{-\alpha}}+$$
$$+\frac{\big(3(3d^2-3c^2-\sqrt3d)C+c(6d+\sqrt3)F\big)\log t}
{18c(3d^2-c^2)rst\sqrt{-\alpha}}=$$
$$=\frac{-(2\sqrt3d+1)\Big(9c^2-\big(\sqrt3d-1\big)^2\Big)}{18c(3d^2-c^2)rst\sqrt{-\alpha}}
\Big(-\frac19P\log r+\frac1{18}(P-Q)\log s+\frac1{18}(P+Q)\log t\Big).$$
Using the first equation in (5.2), we obtain
$$\frac{\partial\vol}{\partial\omega}=\frac{27(-\alpha)}{324c\left(3d^2-c^2\right)rst\sqrt{-\alpha}}
\left(P\log\frac{st}{r^2}+Q\log\frac ts\right)=\frac{\sqrt{-\alpha}}
{12c\left(3d^2-c^2\right)rst}\left(P\log\frac{st}{r^2}+Q\log\frac ts\right).\eqno{\hfill_\blacksquare}$$

\medskip

{\bf6.7.~Proposition.} Let $(\alpha,\omega)\in\overset{\;\circ}\to S$. Then
$\displaystyle\frac{\partial\vol}{\partial\alpha}(\alpha,\omega)<0$ and
$\displaystyle\frac{\partial\vol}{\partial\omega}(\alpha,\omega)>0$.

\medskip

{\bf Proof.} By Proposition 6.6, the inequality
$\displaystyle\frac{\partial\vol}{\partial\alpha}(\alpha,\omega)<0$
follows from
$$\displaystyle M\log\frac{st}{r^2}+N\log\frac{t}{s}<0.\leqno{(6.8)}$$
Indeed, the factor
$\displaystyle\frac{\sqrt3}{648c\left(3d^2-c^2\right)rst\sqrt{-\alpha}}$ in the expression for
$\displaystyle\frac{\partial\vol}{\partial\alpha}(\alpha,\omega)$ is positive because
$(\alpha,\omega)\in\overset{\;\circ}\to S$ implies $c>0$ and $3d^2-c^2>0$ (as observed in the proof of
Lemma 6.3).

Let us establish inequality (6.8). A straightforward calculation shows that
$$M=-54\sqrt3r(t-s)\big(r^2+s^2+t^2-2r(s+t)\big),$$
$$N=-54\sqrt3\big(-2s(s-t)^2t+r^3(s+t)-2r^2(s^2+st+t^2)+r(s^3+s^2t+st^2+t^3)\big)$$
(to verify the above equalities, one may substitute $r,s,t$ respectively by
$\frac{1-\sqrt{3}d}{3},\frac{2-3c+\sqrt{3}d}{6},\frac{2+3c+\sqrt{3}d}{6}$).
Since $r,s,t$ satisfy the triangle inequalities and $0<r<s<t$ due to $(r,s,t)\in\overset\circ\to\Delta_1$,
we have
$$a:=\frac sr>1,\quad0<b:=\displaystyle\frac{t-s}r<1.$$
Hence,
$$M=-54\sqrt3r^4b(2a^2+2ab-4a+b^2-2b+1)$$
$$N=-54\sqrt3r^4(4a^3-2a^2b^2+6a^2b-6a^2-2ab^3+4ab^2-6ab+2a+b^3-2b^2+b)$$
(to verify the above equalities, it suffices to substitute $a,b$ by the respective expressions in $r,s,t$).
It~follows that
$$M\log\frac{st}{r^2}+N\log\frac ts=(M+N)\log\frac tr+(M-N)\log\frac sr=
(M+N)\log(a+b)+(M-N)\log a=-54\sqrt3r^4k$$
where
$$k:=2\big((a+b)(a-1)(2a-b^2+2b-1)\log{(a+b)}-a(a+b-1)(2a-b^2-1)\log a\big)=$$
$$=4b(a+b)(a-1)\log(a+b)+2(2a-b^2-1)\big((a+b)(a-1)\log{(a+b)}-a(a+b-1)\log a\big).$$
Since $4b(a+b)(a-1)\log(a+b)>0$ and $2a-b^2-1>0$ (due to $a>1$ and $0<b<1$), it remains to show that
$$k':=(a+b)(a-1)\log{(a+b)}-a(a+b-1)\log a>0\leqno{(6.9)}$$ Using the series
$$\log x=\sum_{n=1}^{\infty}{\frac{1}{n}\left(\frac{x-1}{x}\right)^n}$$
which converges (absolutely) for each $x\geqslant1$, we have
$$k'=(a+b)(a-1)\sum_{n=1}^{\infty}{\frac1n\frac{(a+b-1)^n}{(a+b)^n}}-a(a+b-1)\sum_{n=1}^{\infty}
{\frac1n\frac{(a-1)^n}{a^n}}=$$
$$=(a-1)(a+b-1)\sum_{n=1}^{\infty}{\frac1n\frac{(a+b-1)^{n-1}}{(a+b)^{n-1}}}-(a-1)(a+b-1)
\sum_{n=1}^{\infty}{\frac1n\frac{(a-1)^{n-1}}{a^{n-1}}}=$$
$$=(a-1)(a+b-1)\sum_{n=1}^{\infty}{\frac1n\left(\left(\frac{a+b-1}{a+b}\right)^{n-1}-
\left(\frac{a-1}a\right)^{n-1}\right)}>0$$
because $\displaystyle\frac{a+b-1}{a+b}-\frac{a-1}{a}=\frac b{a(a+b)}>0$ and
$\displaystyle\frac{a+b-1}{a+b},\frac{a-1}a>0$.

Concerning the derivative with respect to $\omega$, we need to show that
$$P\log\frac{st}{r^2}+Q\log\frac ts>0.$$
The procedure is analogous to the one we just used. First, one can readily verify that
$$P=-6r(t-s)=-6r^2b$$
$$Q=6(2st-rs-rt)=6r^2(2a^2+2ab-2a-b).$$
It follows that
$$P\log\frac{st}{r^2}+Q\log\frac ts=6r^2\Big(-b\log a(a+b)+(2a^2+2ab-2a-b)\log\frac{a+b}a\Big).$$
Finally,
$$-b\log a(a+b)+(2a^2+2ab-2a-b)\log\frac{a+b}a=2k'>0,$$
where $k'$ is defined in (6.9).
$\hfill_\blacksquare$

\medskip

{\bf6.10.~Theorem [\rm{Seidel's Speculation 4}].}~{\it The volume\/ $\vol:S\to\Bbb R$,
$\vol=\vol(\alpha,\omega)$, is decreasing\footnote{Or, as Seidel states, increasing in $|\alpha|=-\alpha$.}
in\/~$\alpha$ and increasing in\/ $\omega$ when\/ $\alpha\ne0$. {\rm(}Tetrahedra corresponding to\/
$\alpha=0$ are degenerate and have vanishing volume{\rm.)}}

\medskip

{\bf Proof.} In view of Proposition 6.7, it suffices to show that the intersection of $S$ with a
horizontal (respectively, vertical) line is, if non-empty, a closed interval.

Remark 3.2.8 and Equations (5.2) imply that the hypotenuse
$\Big\{(c,d)\in\widetilde\Delta_1\mid c=0,0\leqslant d\leqslant\frac{\sqrt3}3\Big\}$ of
$\widetilde\Delta_1$ is sent to
$$c_1:=\bigg\{(\alpha,\omega)\in S\mid\alpha=-\frac1{27}\big(2(6\omega-2)^{\frac32}-3(6\omega-2)+1\big),
\frac13\leqslant\omega\leqslant\frac12\bigg\}$$
by the homeomorphism $\widetilde\Delta_1\to S$, $(c,d)\mapsto\big(\det G_T,\sqrt{\per G_T}\big)$ (see
Corollary 4.5). Similarly, the legs
$\Big\{(c,d)\in\widetilde\Delta_1\mid d=\frac{\sqrt3}3c,0\leqslant c\leqslant\frac14\Big\}$ and
$\Big\{(c,d)\in\widetilde\Delta_1\mid d=-\sqrt3c+\frac{\sqrt3}3,0\leqslant c\leqslant\frac14\Big\}$ of
$\widetilde\Delta_1$ are respectively sent to
$$c_2:=\bigg\{(\alpha,\omega)\in S\mid\alpha=-\frac1{27}\big(-2(6\omega-2)^{\frac32}-3(6\omega-2)+1\big),
\frac13\leqslant\omega\leqslant\frac38\bigg\},$$
$$c_3:=\bigg\{(\alpha,\omega)\in S\mid\alpha=0,\frac38\leqslant\omega\leqslant\frac12\bigg\}.$$
It is not difficult to see that $S$ is the closed region in $\Bbb R^2$ bounded by the curves $c_1,c_2,c_3$;
in other words,
$$S=\bigg\{(x,y)\in\Bbb R^2\mid\frac13\leqslant y\leqslant\frac38,f_1(y)\leqslant x\leqslant f_2(y)\bigg\}
\cup\bigg\{(x,y)\in\Bbb R^2\mid\frac38\leqslant y\leqslant\frac12,f_1(y)\leqslant x\leqslant0\bigg\}
$$
where
$$f_1(y):=-\frac1{27}\big(2(6y-2)^{\frac32}-3(6y-2)+1\big),\quad\frac13\leqslant y\leqslant\frac12,$$
$$f_2(y):=-\frac1{27}\big(-2(6y-2)^{\frac32}-3(6y-2)+1\big),\quad\frac13\leqslant y\leqslant\frac38$$
(see the picture at the end of Section 4).
It follows that, for each $y_0\in\big[\frac13,\frac12\big]$, the intersection of the horizontal line $y=y_0$
with $S$ is a closed interval (it is a single point when $y=\frac13$). As can be easily seen by taking
derivatives, the functions $f_1$ and $f_2$ are invertible (in the indicated domains) and, therefore, a
similar argument involving the inverses of $f_1$ and $f_2$ shows that the intersections of vertical lines
with $S$, if non-empty, are closed intervals.
\hfill$_\blacksquare$

\medskip

It should be noted that the original fourth conjecture states that the volume is {\it decreasing\/} in the
permanent while it is actually increasing. However, this is irrelevant; as discussed at the introduction,
the whole point is to express the volume as a {\it monotonic\/} function of algebraic expressions. In what
follows, we will see that, without an explicit description of the space of ideal tetrahedra modulo
isometries in terms of the determinant and (the square root of the) permanent of doubly stochastic Gram
matrices (that is, the region $S$ in the above theorem), one may be easily led to believe that the volume is
decreasing in the permanent.

\medskip

Let $I_n$ denote the $n\times n$ identity matrix and let $J_n$ denote $n\times n$ matrix whose entries are
all equal to $1$. Let $Z_n$ denote the set of all $n\times n$ symmetric doubly stochastic matrices
which have vanishing diagonal. Van der Waerden's conjecture, proved by Egoritsjev and Falikman, states that
the permanent of a doubly stochastic $n\times n$ matrix attains a unique minimum for the matrix $\frac1nJ_n$
(see [VaL]). In view of this, Seidel formulated a similar conjecture in the context of ideal hyperbolic
tetrahedra. More specifically, Seidel's third conjecture states that the permanent {\it in\/} $Z_n$ attains a
unique minimum for the matrix $\frac1{n-1}(J_n-I_n)$. This conjecture is very simple to prove in our case,
i.e., when $n=4$:

\medskip

{\bf6.11.~Theorem {\rm [Seidel's Speculation 3]}.}~{\it The matrix\/ $\frac13(J_4-I_4)$ is the unique matrix
for which the permanent of the matrices from\/ $Z_4$ attains its minimal value.}

\medskip

{\bf Proof.} Let $M\in Z_4$ be a symmetric doubly stochastic matrix with vanishing diagonal, that is,
$M=\left[\smallmatrix
0&r&s&t\\
r&0&t'&s'\\
s&t'&0&r'\\
t&s'&r'&0\\
\endsmallmatrix\right]$ and $r,s,t,r',s',t'$ are non-negative numbers satisfying $r+s+t=1$, $r+s'+t'=1$,
$r'+s+t'=1$, and $r'+s'+t=1$. One easily concludes that $r=r'$, $s=s'$, and $t=t'$. Hence,
$\per M=\left(r^2+s^2+t^2\right)^2$ and $r+s+t=1$ implies that $\per M$ is minimal when $r=s=t=\frac13$.~
\hfill$_\blacksquare$

\medskip

Note that the space $\Delta$ is strictly contained in $Z_4$ (since, in $Z_4$, it is not required that
$r,s,t$ satisfy the triangle inequalities).

\smallskip

Similarly, we have the following

\medskip

{\bf6.12.~Proposition.}~{\it The matrix\/ $\frac13(J_4-I_4)$ is the unique matrix for which the determinant
of the matrices from\/ $Z_4$ attains its minimal value.}

\medskip

{\bf Proof.}~Let $M$ be a matrix in $Z_4$. As in the proof of Theorem 6.11,
$$\det M=\det\left[\matrix
0&r&s&t\\
r&0&t&s\\
s&t&0&r\\
t&s&r&0\\
\endmatrix\right]=-(r+s+t)(-r+s+t)(r-s+t)(r+s-t)$$
where $r,s,t$ are non-negative numbers such that $r+s+t=1$.

If the numbers $r,s,t$ do not satisfy the triangle inequalities, we have one and only one of the following:
$r>s+t$, $s>r+t$, or $t>r+s$. Indeed, assuming (for example) that $r>s+t$ and $s>r+t$, we obtain
$t<0$, a contradiction. Therefore, if $r,s,t$ do not satisfy the triangle inequalities, then
$\det M=-(r+s+t)(-r+s+t)(r-s+t)(r+s-t)\geqslant0$.

Assume that the numbers $r,s,t$ satisfy the triangle inequalities. Let $A$ denote the area of the Euclidean
triangle of sides of lengths $r,s,t$. Then
$$A=\frac14\sqrt{(r+s+t)(-r+s+t)(r-s+t)(r+s-t)}=\frac14\sqrt{-\det M}.$$
So, $\det M=-16A^2\leqslant0$. Since the area of triangles with fixed perimeter $r+s+t=1$ has a unique
maximum at the equilateral triangle of side lengths $r=s=t=\displaystyle\frac13$, the function
$Z_4\to\Bbb R$, $M\mapsto\det M$, has a unique minimum at the matrix
$\displaystyle\frac13(J_4-I_4)$.
\hfill$_\blacksquare$

\medskip

The fact that the ideal tetrahedron of {\it maximal\/} volume corresponds to the point in $S$ where the
permanent is {\it minimal\/} may be a source of the idea that the volume should be decreasing in the
permanent. However, being on the point $(-\frac1{27},\frac13)\in S$ which corresponds to the ideal
tetrahedron of maximal volume it is not possible to vary (inside of $S$) the permanent while keeping the
determinant constant (see the picture at the end of Section 4).

At the end of the paper the reader may find a couple of graphs of the volume as a function of the
determinant and permanent.

\bigskip

\centerline{\bf7.~Schur functors}

\medskip

We follow Weyl's construction of Schur functors as presented in [Ful] (including the notation).

\smallskip

Let $\bold{FinLin}$ stand for the category of finite dimensional real (or complex) linear spaces. Let
$n\in\Bbb N$ be a natural number and let $\lambda=(\lambda_1,\dots,\lambda_k)$,
$\lambda_1\geqslant\dots\geqslant\lambda_k\geqslant1$, $\lambda_1+\dots+\lambda_k=n$, be a partition of~$n$.
Taking a Young tableau (say, the canonical one) related to $\lambda$, we obtain the corresponding Young
symmetrizer $c_\lambda\in A$, where $A$ denotes the group algebra of the symmetric $n$-group $S_n$. Given a
finite dimensional linear vector space $V$, the symmetric group acts on the right on the tensor power
$V^{\otimes_n}$ by permuting factors. The image of the Young symmetrizer $c_\lambda$ on $V^{\otimes_n}$ is a
linear subspace $\Bbb S_\lambda V\leqslant V^{\otimes_n}$.

The linear space $\Bbb S_\lambda V$ is made up of finite linear combinations of {\it indecomposable\/}
elements of the form
$$v_1\dots v_n:=(v_1\otimes\dots\otimes v_n)c_\lambda.$$
For instance, in the case of the partition $\lambda=(1,\dots,1)$, one gets the exterior power
$\Bbb S_\lambda V=\bigwedge^n V$ whose indecomposable elements are typically denoted by
$v_1\wedge\dots\wedge v_n:=v_1\dots v_n$. In the case of the partition $\lambda=(n)$, the symmetric
power $\Bbb S_\lambda V=S^nV$ is obtained.

We have just arrived at the Schur functor $\Bbb S_\lambda:\bold{FinLin}\to\bold{FinLin}$. At the level of
objects, $V\mapsto\Bbb S_\lambda V$; at the level of morphisms, given a linear map $f:V\to W$ between finite
dimensional linear spaces, then $\Bbb S_\lambda f:\Bbb S_\lambda V\to\Bbb S_\lambda W$ is defined, in terms
of indecomposable elements, by $f(v_1\dots v_n):=f(v_1)\dots f(v_n)$.

\medskip

It is curious to observe that, if a finite dimensional linear space $V$ is equipped with a bilinear
symmetric form $\langle-,-\rangle$ (or a hermitian form in the complex case), then there is an induced form
on $\Bbb S_\lambda V$. It is defined, in terms of indecomposable elements, by
$$\langle v_1\dots v_n,v_1'\dots v_n'\rangle:=\sum_{\sigma\in S_n}\chi(\sigma)g_{1\sigma(1)}g_{2\sigma(2)}
\dots g_{n\sigma(n)},$$
where $g_{ij}:=\langle v_i,v_j'\rangle$ and $\chi$ is the character of the representation of $S_n$ on its
group algebra $A$ induced by the partition $\lambda$. In the exterior power case, this induced form is
nothing but the determinant of the matrix $G:=[g_{ij}]$; in the symmetric power case, it is the permanent of
$G$. In general, $\langle v_1\dots v_n,v_1'\dots v_n'\rangle$ is called the {\it immanant\/} of $G$
(immanants of matrices were introduced in [LiR]).

\smallskip

It is quite common to find connections between the exterior power functor and hyperbolic geometry. Such
connections can be seen, for instance, in Section 2.2 or in [AGr3]. Perhaps, Seidel's
conjectures provide, via the permanent of doubly stochastic Gram matrices of labelled ideal tetrahedra, the
first link of the symmetric power functor with hyperbolic geometry. It is not unreasonable to expect the
other immanants to play a role in hyperbolic geometry, say, in generalizations of Seidel's conjectures to
other (not necessarily ideal) polyhedra or to higher dimensions.

\bigskip

\centerline{\bf8.~A couple of graphs}

\bigskip

\noindent
$\epsfxsize=7cm\vcenter{\hbox{\epsfbox{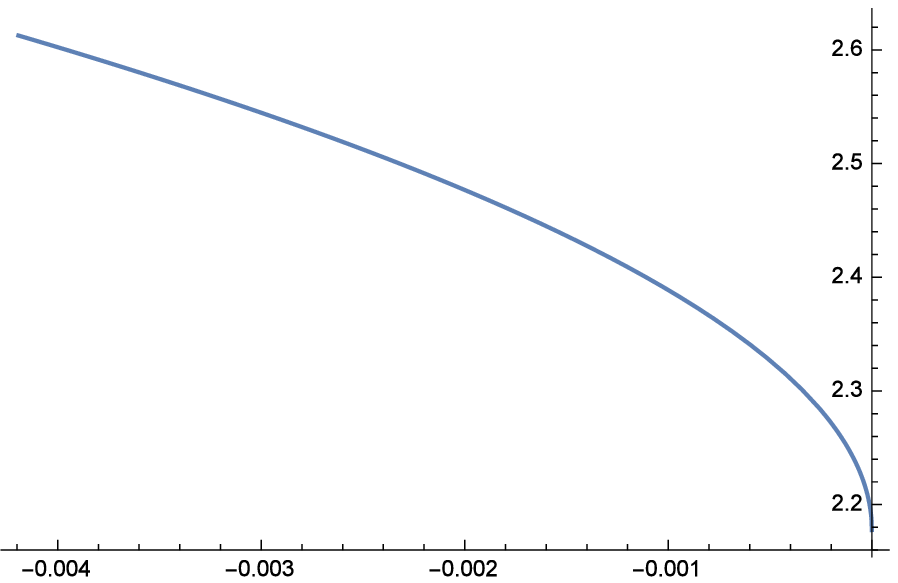}}}$

\medskip

\rightskip255pt

\noindent
Volume as a function of the determinant. The permanent is constant and its square root equals $\frac7{16}$.
The determinant varies between $\frac{14-5\sqrt{10}}{432}$ (isosceles tetrahedron) and $0$ (degenerate
tetrahedron).

\vskip-195.5pt\leftskip250pt

\noindent
$\epsfxsize=7cm\vcenter{\hbox{\epsfbox{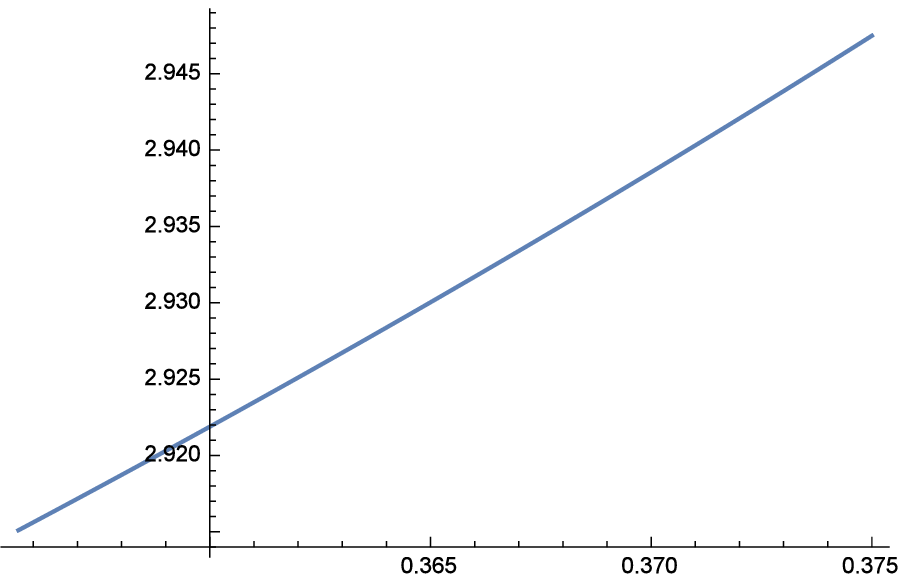}}}$

\rightskip0pt\vskip6pt

\noindent
Volume as a function of the square root of the permanent. The determinant is constant and equals
$-\frac1{54}$. The square root of the permanent varies between $\frac{6-\sqrt3}{12}$ (isosceles
tetrahedron) and $\frac38$ (isosceles tetrahedron).

\leftskip0pt\rightskip0pt

\centerline{\bf References}

\medskip

[Abr] N.~V.~Abrosimov, {\it Seidel's problem on the volume of a non-euclidean tetrahedron,} Doklady
Mathematics, 2010 {\bf82}, 843--846

\smallskip

[AGG] S.~Anan$'$in, C.~H.~Grossi, N.~Gusevskii, {\it Complex hyperbolic structures on disc bundles over
surfaces,} International Mathematics Research Notices, 2011 {\bf19}, 4295--4375

\smallskip

[AGr1] S.~Anan$'$in, C.~H.~Grossi, {\it Basic coordinate-free non-euclidean geometry,} draft of a book
(2011), available at https://arxiv.org/abs/1107.0346

\smallskip

[AGr2] S.~Anan$'$in, C.~H.~Grossi, {\it Coordinate-free classic geometries,} Moscow Mathematical Journal,
2011 {\bf11}, 633--655

\smallskip

[AGr3] S.~Anan$'$in, C.~H.~Grossi, {\it Differential geometry of grassmannians and Plucker map,} Central
European Journal of Mathematics, 2012 {\bf3}, 873--884

\smallskip

[Ana1] S.~Anan$'$in, {\it Complex hyperbolic equidistant loci,} arXiv:1406.5985

\smallskip

[Ana2] S.~Anan$'$in, {\it Reflections, bendings, and pentagons,} arXiv:1201.1582

\smallskip

[Ful] W.~Fulton, J.~Harris, {\it Representation theory, a first course,} Springer-Verlag, 1991

\smallskip

[Huy] D.~Huybrechts, {\it Complex geometry, an introduction,} Universitext, Springer-Verlag, 2004

\smallskip

[LiR] D.~E.~Littlewood, A.~R.~Richardson, {\it Group characters and algebra,} Philosophical Transactions of
the Royal Society A, 1934 {\bf233}, 99--141

\smallskip

[Mil] J.~Milnor, {\it Hyperbolic geometry: the first 150 years,} Bulletin of the American Mathematical
Society, 1982 {\bf6}, 9--24

\smallskip

[Moh] Y.~Z.~Mohanty, {\it Hyperbolic polyhedra: volume and scissors congruence,} University of California,
San Diego, 2002 (Ph.~D. Thesis)

\smallskip

[Sei] J.~J.~Seidel, {\it On the volume of a hyperbolic simplex,} Studia Scientiarum Mathematicarum
Hungarica, 1986 {\bf21}, 243--249

\smallskip

[VaL] J.~H.~Van~Lint, {\it The Van der Waerden conjecture: two proofs in one year,} The Mathematical
Intelligencer, 1982 {\bf4}, 72--77

\enddocument